\input amstex
\documentstyle{amsppt}
\catcode34=12 
\magnification 1200 \TagsOnRight \pageheight{24truecm}
\pagewidth{17truecm} \pageno1 \nopagenumbers

\def\myfootline{\hss\tenrm\folio}
\headline={}
\footline={\ifnum\pageno>1\myfootline\fi} \NoRunningHeads

\expandafter\redefine\csname logo\string @\endcsname{}

\textfont2=\tensy \scriptfont2=\sevensy \scriptscriptfont2=\fivesy
\def\cal{\fam2}
\def\eqalign{\null\vbox\bgroup\advance\baselineskip3pt\halign
\bgroup\hfil${}##{}$&&${}##{}$\hfil\crcr}
\def\endeqalign{\crcr\egroup\egroup\,}

\vskip1cm \topmatter
\def\hat{\widehat}
\newcount\chapno    
\newcount\claimno   
\newcount\stateno   
 \chapno=0
 \stateno=0
 \claimno=0
\define\chap{\global\claimno=0\global\stateno=0\global\advance\chapno
by 1
\ \the\chapno.\enskip}

\define\Clabel{{\global\advance\claimno by 1}\number\chapno.\number\claimno}
\define\clabel#1{\global\advance\claimno by 1\expandafter\xdef
 \csname Cl#1\endcsname{\the\chapno.\the\claimno}\the\chapno.\the\claimno}

\define\Elabel{\tag\global\advance\stateno by 1\the\chapno.\the\stateno}
\define\elabel#1{\tag\global\advance\stateno by 1\expandafter\xdef
 \csname St#1\endcsname{\thetag
       {\the\chapno.\the\stateno}}\the\chapno.\the\stateno}

 \define\cref#1{\expandafter\csname Cl#1\endcsname}
 \define\eref#1{\expandafter\csname St#1\endcsname}

\define\mps#1#2#3{#1\colon #2\to #3}
\define\Mps#1#2#3{#1\colon #2\hookrightarrow #3}
\define\pam#1#2#3#4{#1\hookleftarrow #2@>#3>>#4}

\def\>{\prec}
\def\<{\succ}

\def\dim{\op{dim}}
\def\diam{\op{diam}}

\def\C{{\cal C}}

\def\Šk{{\cal C}}

\def\L{{\cal L}}

\def\P{{\cal P}}
\def\Q{{\cal Q}}

\def\U{{\cal U}}
\def\V{{\cal V}}

\def\W{{\cal W}}

\def\nw#1{{\it {#1}}}
\define\enddem{$\square$\enddemo}
\define\op#1{\operatorname {#1}}
\def\dist#1#2{\op{dist}(#1,#2)}
\def\cov{\op{cov}}

\def\dim{\op{dim}}

\def\Int{\op{Int}}

\def\mean{\rightleftharpoons}

\def\tilde{\widetilde}
\def\phi{\varphi}
\def\matop#1#2#3{\mathop{#1}\limits_{#2}^{#3}}

\def\ii{^{-1}}

\def\Id{\op{Id}}

\def\ANE{\op{ANE}}
\def\AE{\op{AE}}


\topmatter

\abstract
We prove that Dranishnikov's $k$-dimensional
resolution $d_k\colon \mu^k\to Q$ is a \ UV$^{n-1}$-divider  of
Chigogidze's $k$-dimensional resolution $c_k$. This fact implies
that $d_k^{-1}$ preserves $Z$-sets. A further development of the
concept of \ UV$^{n-1}$-dividers permits us to find sufficient
conditions for $d_k^{-1}(A)$ to be homeomorphic to the
N\"{o}beling space $\nu^k$ or the universal pseudoboundary
$\sigma^k$. We also obtain some other applications.
\endabstract

\title
 Preserving  $Z$-sets by Dranishnikov's resolution  \endtitle 
 
\author S. M.  Ageev, M. Cencelj and D. Repov\v s \endauthor
 
\address
Faculty of Mechanics and Mathematics, Belorus State University, 
Av. F. Skorina 4, Minsk, Belarus 220050.\endaddress
\email ageev.sergei\@yahoo.com \endemail

\address
Institute for Mathematics, Physics and Mechanics, and
Faculty of Education,
University of Ljubljana, P.O.B. 2964, Ljubljana,
Slovenia 1001.\endaddress
\email matija.cencelj\@guest.arnes.si \endemail

\address
Institute for Mathematics, Physics and Mechanics, and
Faculty of Education,
University of Ljubljana, P.O.B. 2964, Ljubljana,
Slovenia 1001.\endaddress
\email dusan.repovs\@guest.arnes.si \endemail



\thanks The first author was supported by a grant from
 the Ministry of Education of the Republic of Belarus.
 The second and third author were supported by a grant from
 the Slovenian Research Agency.
 \endthanks

\endtopmatter

\document

\head\chap Introduction
\endhead

 The construction by Dranishnikov \cite{Dr} of the $(k-1)$-soft map
$\mps {d_k}{\mu^k}Q$ of the $k$-dimensional Menger compactum onto
the Hilbert cube $Q$ is an important technique of geometric
topology. This map (called in the sequel the {\it resolution})
permits us to demonstrate the wide analogy between Menger theory
and $Q$-manifold theory. For example,  with its help the
Triangulation and Stability Theorems in the Menger manifold theory
were formulated and proved \cite{Dr}. On the other hand,
Dranishnikov's resolution is by its properties the
finite-dimensional analogue of the projection $\mps{d}{Q\cong
Q\times Q}{Q}$ of the product onto factor with the exception of
being $k$-soft. It is clear that the more properties of
Dranishnikov's resolution will be found the more convenient
instrument it will become.

Chigogidze's resolution  $\mps {c_k}{\nu^k}Q$, which is a bridge
between the N\"{o}beling and Hilbert ($l_2$-) manifold theories,
possesses properties which are every bit as remarkable as those of
Dranishnikov's resolution  (the precise definitions of these
resolutions are put in the Epilogue). The opinion to consider
$c_k$ as the finite-dimensional analogue of the projection
$\mps{c}{\op{l}_2\cong Q\times (-1,1)^\omega}{Q}$  is justified to
a greater degree than in the case of Dranishnikov's resolution.

In \cite{ARS} the investigation of interconnection between these
two resolutions was begun, and it was established that
Chigogidze's resolution is densely contained in Dranishnikov's
resolution, i.e. there exists an embedding
$i_k\colon\nu^k\hookrightarrow\mu^k$ such that
$c_k=d_k\restriction_{\nu^k}$ and $\op{Cl}\nu^k=\mu^k$. This
result is in complete accordance
with the infinite-dimensional situation:
$c=d\restriction_{\op{l}_2}$ and $\op{Cl}\op{l}_2=Q$ for the
natural embedding $i\colon\op{l}_2\cong Q\times
(-1,1)^\omega\hookrightarrow Q$. 

However, while
most of the
useful properties of infinite-dimensional objects ($c,d$ or $i$)
are
evident (for instance, that $i$ is a $\op{UV}$-map), all new
properties of its finite-dimensional analogues are established
with excessive difficulties, especially since $d_k$ fails to be
$k$-soft. The aim of the present
paper is to make a definite progress
towards the investigation of the
finite-dimensional resolutions. Our
main result is:

 \proclaim{Theorem \clabel{1+}} Dranishnikov's resolution $d_k$ is
a {\it $\op{UV}^{k-1}$-divider} of Chigogidze's resolution $c_k$,
i.e. there exists a $\op{UV}^{k-1}$-embedding
$i_k\colon\nu^k\hookrightarrow\mu^k$ such that
$c_k=d_k\restriction_{\nu^k}$ and $\op{Cl}\nu^k=\mu^k$.
 \endproclaim
 The proof of this central theorem is based on a careful analysis
 of the concept of the
 $\op{UV}^{k-1}$-dividers, which may in fact, be
considered as the other
purpose of this paper. 
In particular,
we find a piecewise linear version of Theorem \cref{1+} which is a
crucial ingredient of its proof.

 \proclaim{Theorem \clabel{1+!}}
For every  polyhedron $P$ with the triangulation $L$, there exist
a polyhedron $D$ and  maps $\mps{p}{D}{P}$ and
$\mps{q}{D}{P^{(k)}}$ such that
\item\item{$(1)$} $p$ is a $k$-conservatively soft
map;
\item\item{$(2)$} $p$ is a $\op{UV}^{k-1}$-divider of
a $k$-soft map; and
\item\item{$(3)$} $q$ is
an $\omega$-map for $\omega\>f\ii(L\circ L)$.
 \endproclaim
 
 Theorem \cref{1+} implies several important results. Since the
passing to the preimage with respect to an $n$-soft map preserves
$Z_n$-sets, as it does also with respect to a $\op{UV}^{k-1}$-divider of Chigogidze's
resolution $c_k$ (see Proposition 3.4), Dranishnikov's resolution
{\it preserves $Z$-sets}.

 \proclaim{Theorem \clabel{Def--1}} For
every $F\subset_{Z_k}Q$, $(d_k)\ii(F)\subset_Z\mu^k$.
 \endproclaim
 
 Of course, this intriguing fact will play an important role
in the theory of Menger manifolds. It should also
be remarked  that this fact was
earlier  erroneously claimed to be false \cite{C}.

 Next, we find sufficient conditions for the preimage $d_k\ii(Z)$ of
$Z\subset Q$ to be homeomorphic to the N\"{o}beling space $\nu^k$
or universal pseudoboundary $\sigma^k$.

 \proclaim{Theorem \clabel{1}}
If a complete subset $Z\in\AE(k)$ is strong  $k$-universal with
respect to Polish spaces and $Z\hookrightarrow Q\in\op{UV}^{k-1}$,
then $d_k\ii(Z)\cong\nu^k$.
 \endproclaim
 
 \proclaim{Theorem \clabel{2}}
If a $\sigma$-compact subset $Z\in\AE(k)$ is discretely
$I^k$-approximated and strongly  $k$-universal with respect to
compacta, and also  $Z\hookrightarrow Q\in\op{UV}^{k-1}$, then
$d_k\ii(Z)\cong\sigma^k$.
\endproclaim

 Theorem \cref{1} implies that
$d_k\ii( (-1,1)^\omega)\cong \nu^k$ which is the affirmative
solution of Problem 613 from \cite{MR}. Another
applications of the derived theory will be discussed in the
epilogue.

 \head\chap Preliminaries
\endhead
 Throughout this paper we shall assume all spaces to be
{\it separable complete metric} and all maps to be {\it continuous},
if they do not arise as a result of some constructions and their
properties should be established in the process of proof. The set
of all open covers of the space  $X$ is denoted by $\op {cov}X$.
We will use $\op N(A;\omega )$ to denote the
{\it neighborhood}
$\cup \{U| U\in \omega , U\cap A \ne \emptyset \}$ of $A\subset X$
with respect to $\omega \in \op {cov}X$. The  {\it body } $\cup
\omega$ of a system of sets $\omega$ is the set $\bigcup \{U|
U\in\omega \}$ (we use the sign  $\mean$ for introducing new
objects to the left of it). We say that the embedding $A\subset B$
is {\it strong} and write $A\Subset B$ if $\op{Cl}A\subset\Int B$.

 The refinement of the cover $\omega$ in $\omega _1$
is denoted by $\omega\>\omega _1$. If $\mps {f,g}XY$ are maps, and
$\delta$ is a family of subsets of $Y$, then  {\it the
$\delta$-closeness of $f$ to $g$} (briefly, $\op{dist}(f,g)\>
\delta$ or $f\matop{\sim}{}{\delta}g$) means that if
$f(x)\not=g(x)$, then  $\{f(x),g(x)\}\subset W\in\delta$. The
restriction of a map $f$ onto a subset  $A$ is denoted by
$f\restriction_{A}$ or simply $f\restriction$ if there is no
ambiguity about $A$.  Since  $f$ is an extension of
$g=f\restriction_{A}$, we write this as $f=\op{ext}(g)$.
If $\delta>0$ is a number, then  $\delta$-closeness of $f$ to $g$
is denoted by $\op{dist}(f,g)<\delta$, as in the case of covers.
We denote the distance $d(x,y)$ between points $x,y\in X$ of
metric space $(X,d)$ as $\dist{x}{y}$ if there is no confusion.

 Let us introduce a series of notions concerned with the extension
{\it of partial maps}, i.e. maps given on closed subspaces of a
metric space \cite{Hu}. If an arbitrary partial map
$\pam{Z}{A}{\phi}{X},\dim Z\le k, k\le\infty$, can be extended on
the entire space $Z$ [on some neighborhood of $A$], then $X$ is
called an {\it absolute [neighborhood] extensor} in dimension $k$,
$X\in\op{A[N]E}(k)$. If $k=\infty$, then the notion of absolute
[neighborhood] extensor ($X\in \op{A[N]E}$) arises. By the
Kuratowski-Dugundji Theorem \cite{Hu}, the property of
extendability in finite dimension correlates with the connectivity
and the local connectivity of the space:
$X\in\op{AE}(k)\Leftrightarrow X\in\op{C}^{k-1}\&\op{LC}^{k-1}$.

 The problem of extension of partially defined morphisms
has a categorical character. In the category of maps having a
fixed target $Y$ the problem of extension of morphisms is known as
the
problem of extension of a partial lift to the global lift. For a
given map $\mps fXY$, the  {\it partial lift of the map } $\mps
{\psi }ZY$ with respect to $f$ is the map $\mps {\varphi }AX$
which is defined on the closed subset $A\subset Z$ and which makes
the following diagram commutative.
 $$ \CD A   @>\varphi >>  X \\ @ViVV
@VVfV \\ Z  @>>\psi >  Y
\endCD
$$
 A partial lift $\varphi$ of the map $\psi$ is
extended to a {\it global (local) lift with respect to $f$} if
there exists a global (local) extension of $\mps {\varphi}ZX$
which is the lift of $\psi$. Thus, the problem of global lifting
consists in the splitting of the square diagram above by the map
$\hat \varphi $ into two triangular commutative diagrams.

Recall that  a map $f$ is called {\it  soft (local soft)} with
respect to pair  $(Z,A)$, if any partial lift $\mps {\phi}AX$
(with respect to $f$) of any map $\mps {\psi}ZY$ can be extended
to the global (local) lift. The collection $\goth S(f)$ of all
pairs $(Z,A)$ for which $f$ is soft  will be called a {\it
softness envelope } of the map $f$. Note that if $|Y|=1$, then
the problem of extension of lifts is transformed into the problem
of extension of maps.

Let $\goth C$ be a class  of pairs  $(Z,A)$   in  which $A\subset
Z$ is a closed subset. The map $f$ is called  {\it $\goth C$-soft
(local $\goth C$-soft)} if it is soft (local soft ) with respect
to  all pairs  $(Z,A)$  from $\goth C$. Along this line, we can
introduce the notions of $(n,k)$-softness  ($\goth C=\lbrace
(Z,A)\mid \op {dim}Z\le n,\op {dim} A\le k\rbrace$), polyhedral
softness ($\goth C=\lbrace (Z,A)\mid Z,A$ are polyhedra, and $\op
{dim}Z\le n\rbrace$). We denote the class of $(n,n)$-soft maps, or
briefly, $n$-soft maps by $\goth S_n$. If $\goth S(f)$ contains
all pairs $(Z\times \lbrack 0,1\rbrack,Z\times \lbrace 0\rbrace)$,
where $\dim Z\le n$, then $f$ is called a {\it Hurewicz
$n$-fibration}.
 The  following assertion is well-known \cite{Hu}:
 
 \proclaim{Proposition \clabel{Prel5+}} Let $Y\in \op {ANE}(n)$.
Then for every  $\varepsilon\in\op{cov}Y$ there exists
$\delta\in\op{cov}Y,\delta\>\varepsilon$ such that  for every
closed subspace $A\subset W,\dim W\le n$, and also  for all maps
$\mps{\hat\alpha}{W}{Y}$ and $\mps{\beta}{A}{Y}$ such that
 $\op{dist}(\hat\alpha \restriction_{A},\beta)\>\delta$, there
exists an extension  $\mps{\hat\beta}{W}{Y}$, $\hat\beta
\restriction_{A}=\beta$, such that
$\op{dist}(\hat\alpha,\hat\beta)\>\varepsilon$.
\endproclaim

   We say that the  family $\cal L$ of closed subsets in metric
space $Y$ is an {\it $\op{equi-LC}^{n-1}$-family} provided that
for any $x\in \cup \cal L$ and $\varepsilon>0$ there exists
$\delta>0$ such that  any map $\mps
{\phi}{S^k}{\op{N}(x;\delta)\cap L}$, $L\in\cal L$ and $k<n$,
defined on the boundary of the ball $B^{k+1}$, is extended to the
map $\mps {\hat \phi}{B^{k+1}}{\op{N}(x;\varepsilon)\cap L}$. The
local $n$-softness of the map  $f\colon X \to Y$ of a complete
space $X$, by the Michael Selection Theorem \cite{SB}, is
equivalent to $\lbrace f^{-1}(y)\rbrace \in \op {equi-LC}^{n-1}$.
The following assertion is a far reaching generalization of the
Michael Selection Theorem.

  \proclaim{Theorem  \clabel{3.3.} (Filtered finite-dimensional
selection theorem \cite{SB})}
 Let $Y$ be a complete metric space, and $X$ a  $k$-dimensional
paracompact space. Let $\L_{i},0\le i\le k$, be an
$\op{equi}-\op{LC}^{k-1}$-family of subsets of $Y$ such that
$\L_{0}\subset \L_{1}\subset \dots\subset\L_{k}$, $F_0\subset
F_1\subset\dots\subset F_{k}$ an increasing sequence of lower
semicontinuous  mappings $F_{k}:X\rightsquigarrow Y$ with closed
values. If $\{ F_i(x)|x\in X\}\subset\L_{i}$ for every $0\le i\le
k$, and the embedding $F_i(x)\hookrightarrow F_{i+1}(x)$ is {\it
$(k-1)$-aspherical} for every  $x\in X$ and for every $0\le i\le
k$ (i.e. this embedding induces trivial homomorphism of homotopy
groups $\pi_i,i\le k-1$), then there exists a selection
$\mps{s}{X}{Y}$ of $F_{k}$.
\endproclaim

 If all $F_i$ are equal to each other, this theorem
is transformed to the $k$-dimensional Michael Selection Theorem.
Studying the interconnection of discrete $I^n$-approximation and
$\op{UV}^n$-division, we need the following easy consequence of
Theorem \cref{3.3.}, the direct proof of which is unknown to us.

 \proclaim{Proposition \clabel{Prel5+-}} Let
$\mps{f}{Y}{Z}$ be an $n$-soft map defined on a complete
$\op{ANE}(n)$-space $Y$. Then for every $\varepsilon\in\op{cov}Y$
there exists $\delta\in\op{cov}Z$ such that for every  space
$W,\dim W\le n$, for all maps $\mps{\alpha}{W}{Y}$ and
$\mps{\beta}{W}{Z}$ with $\op{dist}(f\circ\alpha,\beta)\>\delta$,
there exists a map $\mps{\tilde{\beta}}{W}{Y}$ such that
$f\circ\tilde{\beta}=\beta$ and
$\op{dist}(\alpha,\tilde{\beta})\>\varepsilon$.
\endproclaim

 We say that  a {\it dense map} $\mps{f}{X}{Y}$ (i.e.
$f(X)$ is dense in $Y$)  from  $\op{ANE}(k)$-space $X$ into $Y$
is:
 \item\item{$(1)$} A {\it $\op{UV}^{k-1}$-map}
(briefly,  $f\in\op{UV}^{k-1}$) if for every neighborhood
$\U(y),y\in Y$, there exists a neighborhood $\V(y)$ such that the
embedding $f^{-1}(\V(y))\hookrightarrow f^{-1}(\U(y))$ is
$(k-1)$-aspherical; and
  \item\item{$(2)$} The map $f$ is
{\it approximately polyhedrally $k$-soft} if for every
$\varepsilon>0$ there exists $\delta>0$ such that  for every
$k$-dimensional compact polyhedral pair $(W,A)$ and for all maps
$\mps{\varphi}{A}{X}$ and $\mps{\psi}{W}{Y}$ with
$\dist{\psi}{f\circ\varphi}<\delta$, there exists an extension
$\mps{\hat \varphi}{W}{X}$ of $\varphi$ satisfying
$\op{dist}(\psi,f\circ \hat \varphi)<\varepsilon$.

 In general, a $\op{UV}^{k-1}$-preimage  of $\ANE(k)$-space is not
an $\ANE(k)$-space. But there exists one important exception
\item\item{$(3)$} If $X_0\hookrightarrow X\in\op{UV}^{k-1}$
and $X\in\ANE(k)$, then  $X_0\in\ANE(k)$.

\noindent The  following criterion is well-known (see
\cite{Sch}):
 \proclaim{Proposition \clabel{s-Z3}} If $f\in\op{UV}^{k-1}$, then
$f$ is an approximately polyhedrally  $k$-soft map. Conversely, if
$f$ is an approximately polyhedrally $k$-soft map and
$Y\in\ANE(k)$, then $f\in\op{UV}^{k-1}$.
\endproclaim
 From Proposition \cref{s-Z3} we deduce several
known properties of $\op{UV}^{k-1}$-maps.
 \proclaim{Proposition  \clabel{s-Z3-5}} If $\mps fXY\in
\op{UV}^{k-1}$ and $Y$ is complete, then  $Y\in\ANE(k)$. If
additionally $X\in\AE(k)$, then  $Y\in\AE(k)$.
 \endproclaim
 \proclaim{Proposition \clabel{s-Z5}}
Let $\mps gXY$ be a $\op{UV}^{k-1}$-map, and let $\mps fXZ$ and
$\mps hYZ$ be maps between  $\op{ANE}(k)$-spaces such that
$f=h\circ g$. Then $f\in\op{UV}^{k-1}$ if and only if
$h\in\op{UV}^{k-1}$.
 \endproclaim
 \proclaim{Proposition  \clabel{s-Z3-7}} If $\mps fXY$ is
a $\op{UV}^{k-1}$-map of $\ANE(k)$-spaces, and
$f\ii(Y_0)\hookrightarrow X\in \op{UV}^{k-1}$, then  $\mps
{f\restriction}{f\ii(Y_0)}{Y_0}\in \op{UV}^{k-1}$.
 \endproclaim
 Recall that the {\it fiberwise product} $W=X_f\times _gZ$ of $X$
and $Z$ with respect to  $\mps fXY$ and $\mps gZY$ is the subset
$\lbrace (x,z)\mid f(x)=g(z)\rbrace \subset X\times Z$. The
projections of $X\times Z$ onto  $Z$ and onto $X$ generate the
maps $\mps {f'}WZ$  and  $\mps {g'}WX$ which are called the {\it
projections parallel $f$ and $g$} respectively. We denote it
$f'\parallel f$ and $g'\parallel g$ for brevity.

 Several properties of maps are inherited  by parallel projections.
For instance, the softness envelope $\goth S(f) $ is contained in
$\goth S(f') $. The following is easily established
 \item\item{$(a)$} If $f$ is
$\ n$-soft  and $g\in\op{UV}^{n-1}$, then   $f'\parallel f$ is
$n$-soft , and $g'\parallel g$ is a $\op{UV}^{n-1}$-map.

 In \cite{ARS} we described the sufficiently represented part
of softness envelope  of  Dranishnikov's resolution $d_n$.
  \definition{Definition  \clabel{Prel6}}
The pair $(Z,A)$  is called {\it $n$-conservative} if any partial
lift $\mps {\phi}{A}{S^n\times S^n}$ of $\mps {\psi}{Z}{S^n}$ with
respect to the projection $\mps {\op {pr_2}}{S^n\times S^n}{S^n}$
of the $n$-spheres product onto the second factor is extended to
the global lift $\mps {\hat \phi}Z{S^n\times S^n}$ of  $\psi$ such
that $(\hat \phi)^{-1}(\op {Diag})\subset A.$
 \enddefinition
 The map $\mps fXY$ which is soft with respect to
all $n$-dimensional $n$-conservative pairs   $(Z,A)$ is called
{\it $n$-conservatively soft}. Dranishnikov's resolution $d_n$ is
$n$-conservatively  soft \cite{ARS}. This, in turn, implies all
known soft properties of $d_n$.
 \definition{Definition  \clabel{new-1}} {\it The Polish
space } $X$ (i.e.  complete and separable) is called \nw{strong
$k$-universal Polish}  (or \nw{strong $k$-universal with respect
to  Polish spaces}) if any map $\mps\phi ZX$ of Polish space
$Z,\dim Z\le k$ is arbitrarily closely approximated by closed
embedding.
 \enddefinition
 \definition{Definition  \clabel{new-2}} Let $\{I^k_i\}$ be a
countable family of $k$-dimensional disks, and $D$ -- their
discrete union  $\coprod\{I^k_i\mid 1\le i<\infty\}$. The space
$X$ is called {\it discretely $I^k$-approximated } if any map
$\mps{\phi}{ D}X$ is arbitrarily closely  approximated by a map
$\mps{\tilde\phi}{D}X$ with discrete $\{\tilde\phi(I^k_i)\mid 1\le
i<\infty\}$.
 \enddefinition
 For Polish $\op{ANE}(k)$-space the property of strong
$k$-universality with respect to  Polish spaces and discrete
$I^k$-approximation are equivalent \cite{Bw, p. 127}. In
\cite{Ag,AgM} the following criterion for the N\"{o}beling space
was established.
  \proclaim{Theorem \clabel{sug1}} The Polish space
$X$ of dimension $k$ is homeomorphic to the N\"{o}beling space
$\nu^k$ if and only if $X$ is an $\op{AE}(k)$-space which is
strongly $k$-universal with respect to Polish spaces.
 \endproclaim
 Let ${\cal C}$ be a class  of spaces.
Recall \cite{BM} that  a  space $X$ is called {\it strongly ${\cal
C}$-universal} if any map $\mps{f}{D}{X}$ of $D\in\cal C$, the
restriction of which on a closed subspace $C\in\cal C$ is a
$Z$-embedding, is arbitrarily closely approximable by a
$Z$-embedding $\mps{f'}{D}{X}$ with
$f\restriction_{C}=f'\restriction_{C}$. Under the {\it universal
$k$-dimensional pseudoboundary} we understand a ${\cal
C}_{c(k)}$-absorber  $X$, where ${\cal C}_{c(k)}$ means the class
of all  $k$-dimensional compacta. This means that $k$-dimensional
$\sigma$-compact $\op{AE}(k)$-space $X$ is strongly ${\cal
C}_{c(k)}$-universal and discretely $I^k$-approximated. The paper
\cite{Za} called attention to the fact that Theorem \cref{sug1}
implies the uniqueness of the topological type of the universal
$k$-dimensional pseudoboundary.
 \proclaim{Theorem \clabel{sug1++} (see \cite{CZ} and \cite{Ag})}
Any two universal $k$-dimensional pseudoboundaries are
homeomorphic.
 \endproclaim


\head\chap Basic properties of $\op{UV}^{n-1}$-divider
\endhead
 By $\cal P$ we denote a subclass of  $n$-soft maps of
$\op{ANE}(n)$-spaces.
  \definition{Definition  \clabel{Def1}}
 A proper map $\mps hYZ$ between $\op{ANE}(n)$-spaces
is called
 \item\item{$(i)$} A $\op{UV}^{n-1}$-{\it divider}
of $\mps fXZ$ if there exists a topological embedding $g\colon
X\hookrightarrow Y\in\op{UV}^{n-1}$ such that $f=h\circ g$ (i.e.
$Y_0\mean g(X)$ is dense $G_\delta$ in $Y$, $Y_0\hookrightarrow
Y\in\op{UV}^{n-1}$ and $f=h\restriction_{Y_0}$); and
  \item\item{$(ii)$} A $\op{UV}^{n-1}$-{\it divider}
of  $\cal P$ if $h$ is a $\op{UV}^{n-1}$-divider of some map $\mps
fXZ\in\cal P$.
 \enddefinition
Our interest is  basically in the $\op{UV}^{n-1}$-dividers $\mps
hYZ$ of $\cal P$ with $\dim Y=n\le\dim Z$. The first nontrivial
example of a $\op{UV}^{n-1}$-divider was constructed in \cite{C}.
Prior to establishing that Dranishnikov's resolution is an
$\op{UV}^{n-1}$-divider of Chigogidze's resolution we present
their general properties. It follows from Proposition \cref{s-Z5}
that a $\op{UV}^{n-1}$-divider $h$ of $f$ is $\op{UV}^{n-1}$ iff
$f\in\op{UV}^{n-1}$. Hereafter and also from ${\cal P}\subset
{\goth S}_n$ it easily follows that:
   \proclaim{Proposition \clabel{Def2}}
Any $\op{UV}^{n-1}$-divider of $\cal P$ is an open
$\op{UV}^{n-1}$-map between $\op{ANE}(n)$-spaces.
 \endproclaim
  If in the definition of $\op{UV}^{n-1}$-divider we restricted
ourselves to compact spaces, then the conclusion of Proposition
\cref{Def2} can be essentially strengthened \cite{Sch,ARS}.
 \proclaim{Theorem \clabel{Def3} (on division of locally
$n$-soft maps of compact spaces)} Let the locally $n$-soft map
$\mps fXZ$ be a composition of a $\op{UV}^{n-1}$-map $\mps gXY$
and a map $\mps hYZ$. If all spaces  $X,Y$ and $Z$ are
$\ANE(n)$-compacta, then  $h$ is locally $n$-soft.
 \endproclaim
 By Proposition \cref{s-Z3}, any $\op{UV}^{n-1}$-map is
approximately polyhedrally $k$-soft. On the other hand, the
passing to the preimage with respect to locally $n$-soft map
preserves $Z_n$-sets. These facts easily implies:
 \proclaim{Proposition \clabel{Def-3}} If
$\mps hYZ$ is a $\op{UV}^{n-1}$-divider of  $\cal P$, then
$h\ii(F)\subset_{Z_n}Y$, for every  $F\subset_{Z_n}Z$.
 \endproclaim
 \smallskip
 We say that  the subclass  $\cal P$ is closed  with respect to
 \item\item{$(1)$} {\it Composition} if for any
$\mps{f}{X}{Y}\in\cal P$ and $\mps{g}{Y}{Z}\in\cal P$,
$\mps{g\circ f}{X}{Z}\in\cal P$; and
 \item\item{$(2)$} {\it Passing to complete preimages} if
for any $\mps{f}{X}{Y}\in\cal P$ and  $\ANE(n)$-subspace
$Y_0\hookrightarrow Y$, $\mps{f\restriction_{X_0}}{X_0\mean
f\ii(Y_0)}{Y_0}\in\cal P$.

 \noindent It easily seen that the class  of all  $n$-soft
strongly $n$-universal maps of Polish $\ANE(n)$-spaces (which are
the basic interest of this paper) satisfies the conditions $(1)$
and $(2)$.
 \proclaim{Proposition \clabel{Def4}}
 Let $\cal P$ be a class which is closed both
with respect to  composition and passing to complete preimages. If
$\mps{h_1}{Y_1}{Y_2}$ is a $\op{UV}^{n-1}$-divider of
$\mps{f_1}{X_1}{Y_2}\in\cal P$, where $g_1\colon
X_1\hookrightarrow Y_1\in\op{UV}^{n-1}$, and $\mps{h_2}{Y_2}{Y_3}$
is a $\op{UV}^{n-1}$-divider of $\mps{f_2}{X_2}{Y_3}\in\cal P$,
where $g_2\colon X_2\hookrightarrow Y_2\in\op{UV}^{n-1}$, then the
composition $\mps{h_2\circ h_1}{Y_1}{Y_3}$ is a
$\op{UV}^{n-1}$-divider of the composition $X_0\mean
X_1\cap(h_1)\ii(X_2)\matop{\longrightarrow}{}{f_1\restriction_{X_0}}X_2
\matop{\longrightarrow}{}{f_2}Y_3\in\cal P$.
 \endproclaim
 \demo{Proof of \cref{Def4}} Since
$f_1\restriction_{X_0}\in\cal P$ (as the restriction of $f_1$ onto
the  complete preimage  $X_0=(f_1)\ii(X_2)$) and $f_2\circ
f_1\restriction_{X_0}\in\cal P$ (as the composition of maps from
$\cal P$), it suffices to prove that $e\colon X_0\hookrightarrow
Y_1$ is a $\op{UV}^{n-1}$-embedding. But the embedding $e_1\colon
X_0\hookrightarrow X_1\in\op{UV}^{n-1}$, being a parallel
projection in the fiberwise product of $n$-soft map
$\mps{f_1}{X_1}{Y_2}$ and $\op{UV}^{n-1}$-embedding $g_2\colon
X_2\hookrightarrow Y_2$ (see \cref{s-Z3-7}$(a)$). Then the
embedding $e\in\op{UV}^{n-1}$, being a composition of
$\op{UV}^{n-1}$-maps $e_1$ and $g_1$.
\enddem
  We give one more property of $\op{UV}^{n-1}$-dividers.
 \proclaim{Proposition \clabel{Def6}} Let $\cal P$ be closed
with respect to  passing to complete preimages, and $\mps hYZ$ be
a $\op{UV}^{n-1}$-divider of $\mps{f}{X}{Z}\in\cal P$, where
$g\colon X\hookrightarrow Y\in\op{UV}^{n-1}$. Then for any
$Z_0\hookrightarrow Z\in\op{UV}^{n-1}$ the following holds
 \item\item{$(3)$} $\mps{f\restriction_{X_0}}
{X_0}{Z_0}\in\cal P$ where $X_0\mean f\ii(Z_0)=h\ii(Z_0)\cap X$;
 \item\item{$(4)$} $X_0\hookrightarrow Y\in\op{UV}^{n-1}$;
 \item\item{$(5)$} $\mps{h\restriction_{Y_0}}{Y_0}{Z_0}$ is an
$\op{UV}^{n-1}$-divider of $f$; and
 \item\item{$(6)$} $Y_0\mean h\ii(Z_0)
\hookrightarrow Y\in\op{UV}^{n-1}$.
 \endproclaim
 \demo{Proof of \cref{Def6}} By the $n$-softness of
$f$ it follows  that  $X_0\hookrightarrow X\in\op{UV}^{n-1}$ and
$X_0\in\ANE(n)$. Since $Z\in\ANE(n)$ and $Z_0\hookrightarrow
Z\in\op{UV}^{n-1}$, it follows by \cref{Prel5+-}$(3)$,   that
$Z_0\in\ANE(n)$. Then the conditions imposed on the subclass $\cal
P$ imply that $f\restriction_{X_0}\in{\cal P}$, hence $(3)$ is
proved. Since $X_0\hookrightarrow X$ and $g\colon X\hookrightarrow
Y\in\op{UV}^{n-1}$, Proposition \cref{s-Z5} implies
$X_0\hookrightarrow Y\in\op{UV}^{n-1}$ which proves $(4)$.

 The  property $(5)$ is equivalent  to the  following fact.
 \proclaim{Lemma \clabel{Def31}} $X_0\hookrightarrow
Y_0\in\op{UV}^{n-1}$.
 \endproclaim
 \demo{Proof of \cref{Def31}}
 Consider a neighborhood $\U\subset Y_0$ and a map
of pairs $\mps{\varphi} {(B^n,S^{n-1})}{(\U,\U\cap X_0)}$. By
virtue of $g\in\op{UV}^{n-1}$, the map $\varphi$ can be
arbitrarily closely approximated by $\mps{\varphi'} {B^n}{\U\cap
X}$ with $\varphi'=\varphi$ on $S^{n-1}$. Also, by
$Z_0\hookrightarrow Z\in\op{UV}^{n-1}$, the map  $f\circ \varphi'$
can be arbitrarily closely approximated by  $\mps{\psi}
{B^n}{Z_0}$ with $\psi=f\circ \varphi$ on $S^{n-1}$. And finally,
$n$-softness of $f$ implies the existence of a lifting
$\mps{\tilde\psi} {B^n}{X_0}$ of $\psi$ which coincides with
$\varphi'$ on $S^{n-1}$ and is arbitrarily close to $\varphi'$.
\enddem
It follows by Lemma \cref{Def31} and \cref{s-Z3-5} that
$Y_0\in\op{ANE}(n)$. Since the composition $X_0\hookrightarrow
Y_0\hookrightarrow Y$ is a $\op{UV}^{n-1}$-embedding it follows by
\cref{Def31} and \cref{s-Z5} that $Y_0\hookrightarrow
Y\in\op{UV}^{n-1}$. Hence $(6)$ is proved.
 \enddem
 \smallskip
 Up to the end of the section we shall fix an $n$-dimensional
space  $Y$ and a $\op{UV}^{n-1}$-divider $\mps{h}{Y}Z$ of an
$n$-soft map $\mps{f}{X}{Z}$ (with  $g\colon X\hookrightarrow
Y\in\op{UV}^{n-1}$).
 \proclaim{Proposition \clabel{Def5}}
If $Z$ is discretely $I^n$-approximated, then  $Y$ is also
discretely $I^n$-approxi\-mated.
 \endproclaim
 \demo{Proof } Consider maps $\mps{\varepsilon}{Y}{(0,1)}$
and $\mps{\varphi}{D}{Y}$, where $D$ is a countable discrete union
$\coprod\{I^n_i\mid i<\infty\}$ of $n$-dimensional cubes. Since
$h$ is proper, we can assume that  $\varepsilon$ coincides with
$\zeta\circ f$ for a sufficiently small function
$\mps{\zeta}{Z}{(0,1)}$.

Since $g\in\op{UV}^{n-1}$, $\varphi$ is approximated by a map
$\mps{\varphi'}{D}{X}$ such that
$\dist{\varphi'}{\varphi}\>\varepsilon\circ \varphi$. Next, we
approximate $\psi'\mean f\circ\varphi'$ sufficiently closely by a
map $\mps{\psi}{D}{Z}$ for which the family $\{\psi(I^n_i)\mid
i<\infty\}$ is discrete. By \cref{Prel5+-}, $\psi$ can be lifted
to the map $\mps{\tilde\psi}{D}{X}$ which is arbitrarily close to
$\varphi'$. It can be easily seen that the  family
$\{\tilde{\psi}(I^n_i)\}$ is discrete in $Y$, and $\tilde{\psi}$
is at the required distance from $\varphi$.
 \enddem
 \smallskip
The proof of Theorems \cref{1} and \cref{2} will be given in
Section 4, and the rest of this section presents some necessary
results for this.

Since the notions of strong $n$-universality with respect to
Polish spaces and discrete $I^n$-approximateness  are equivalent
for Polish $\op{ANE}(n)$-spaces, we can assert, using the
criterion of the N\"{o}beling space $\nu^n$ (Theorem \cref{sug1}),
that
 \item\item{$(a)$} If a Polish space $Z\in\AE(n)$ is
strongly $n$-universal with respect to  Polish spaces and $\dim
Y=n$, then  $Y=h\ii(Z)\cong\nu^n$.

  \proclaim{Proposition \clabel{Def55}} Let $Z$ be a
discretely $I^n$-approximated and strong ${\cal
C}_{c(n)}$-universal space, where we denote by ${\cal C}_{c(n)}$
the class  of all $n$-dimensional compacta. Then $Y=h\ii(Z)$ is
strongly ${\cal C}_{c(n)}$-universal.
 \endproclaim
\demo{Proof } By \cref{Def5}, $Y=h\ii(Z)$ is a discretely
$I^n$-approximated space. Since any compactum in a discretely
$I^n$-approximated $\op{ANE}(n)$-space is a $Z_n$-set \cite{BM},
it follows that
 \item\item{$(iii)$} Any compactum  in $Z$ (as well as in $Y$)
is a $Z_n$-set.

Let $\mps{\varphi}{D}{Y}$ be a map of an $n$-dimensional compactum
$D$ such that its restriction onto a closed subspace $C$ is an
embedding. Since $g\in\op{UV}^{n-1}$, we can assume without loss
of generality that $\varphi(D\setminus C)\subset X$. It follows
from $(iii)$ that $(h\circ\varphi)(C)\subset_{Z_n}Z$. Therefore
$h\circ\varphi$ can be arbitrarily closely approximated by a map
$\mps{\psi}{D}{Z}$ such that  $\psi=\varphi$ on $C$, and
$\psi\restriction_{D\setminus C}$ is an embedding whose image does
not intersect $\varphi(C)$.

Let $D_0\mean\varphi\ii(X)\subset D$. Since the map
$\mps{f}{X}{Z}$ is $n$-soft, $\psi\restriction_{D_0}$ can be
lifted with respect to $f$, by Proposition \cref{Prel5+-}, to the
map $\mps{\chi}{D_0}{X}$, arbitrarily close to
$\mps{\varphi\restriction_{D_0}}{D_0}{X}$. It is easy to conclude
that the map $\mps{\varphi'}{D}{Y}$, defined as $\varphi'=\varphi$
on $C$ and  $\varphi'=\chi$ on $D\setminus C$, is continuous. It
is clear that the map $\varphi'$ is an embedding which is
arbitrarily close to $\varphi$. Hence the proof is completed.
 \enddem
 With the help of the criterion of universal
pseudoboundary (Theorem \cref{sug1++}) Proposition \cref{Def55}
implies that
 \item\item{$(b)$} If $Z\in\AE(n)$ is
$\sigma$-compact  discretely $I^n$-approximated and strongly
$n$-universal with respect to  compact spaces, then $Y$ is
homeomorphic to the universal $k$-dimensional pseudoboundary
$\sigma^n$.


\head\chap Inverse limit properties of $\op{UV}^{n-1}$-dividers
\endhead
 \definition{Definition  \clabel{new-3}}
 A map $\mps fXY$ is called {\it strong }  $n$-{\it universal
with respect to  Polish spaces} if for every $n$-dimesional Polish
space $Q$ and also  for every maps $\mps{\varepsilon}{X}{(0,1)}$
and $\mps {\varphi}QX$ there exists a closed embedding
$\varphi'\colon Q\hookrightarrow X$ $\ \varepsilon$-close to
$\varphi$ such that  $f\circ \varphi'=f\circ \varphi$.
 \enddefinition
 \definition{Definition  \clabel{7.6}}
 A map $\mps fXY$ is called
{\it $n$-complete} if for every  map $\mps {\varphi}ZX$ of
$n$-dimensional Polish space $Q$ there exists a closed embedding
$\psi\colon Q\hookrightarrow X$ such that $f\circ \varphi=f\circ
\psi$.
 \enddefinition
 Note that  the composition of  $n$-soft and  $n$-complete
maps is $n$-complete.
 \smallskip
The commutative diagrams  $({\cal D}_t)$ for $t=1,2,3,\dots$
 $$\CD X_{t+1} @>\eta _t>> X_t \\
 @Vf_{t+1}VV  @VVf_tV \\
Z_{t+1} @>\sigma _t>> Z_t
\endCD \tag{${\cal D}_t$} $$
 generate the map $
f\colon X\mean\varprojlim \lbrace X_t,\eta_t\rbrace\rightarrow
Z\mean \varprojlim \lbrace Z_t,\sigma _t\rbrace$ of inverse limit
of spectra. In general, $n$-softness ($n$-conservative softness
and so on) of all maps $f_t$, $\sigma _t$, $\eta_t$  does not
imply that $f$ possesses the corresponding property. As care
should be taken to see that the map properties are preserved by
passage to the inverse limit of spectra, we introduce the
following
 \definition{Definition  \clabel{Prel5}} The commutative
diagram $({\cal D}_t)$ {\it possesses a property $\Q$}, if {\it
its characteristic map } $\mps {\chi_{t}}{X_{t+1}}{W_{t}}$ into
the fiberwise product  $W_{t}\mean (Z_{t+1})_{\sigma _{t}}\times
_{f_{t}} X_{t}$, given by $\chi_{t}(x)=(f_{t+1}(x),\eta_{t}(x))\in
W_{t}$, possesses $\Q$.
\enddefinition
 Basically, we are interested in {\it $n$-soft  and
$n$-complete} commutative diagrams. The particular case of the
following proposition is given in \cite{FC, 2.2.4}.
 \proclaim{Proposition \clabel{7.5.1}} Let
$\mps f{X}{Z}$ be a map of inverse limit of spectra $X\mean
\varprojlim \lbrace Y_t,\eta_t\rbrace$ and $Z\mean \varprojlim
\lbrace Z_t,\sigma _t\rbrace$, generated by commutative diagrams
$({\cal D}_t)$, $t\ge 1$. Let also the diagrams $({\cal D}_t),t\ge
1$, be $n$-soft and $n$-complete, and $f_{1}$ and maps $\sigma
_{t},t\ge 1$, $n$-soft. Then the map $f$ is $n$-soft and strong
$n$-universal with respect to Polish spaces.
\endproclaim
 \demo{Proof of \cref{7.5.1}}
Since  $f$ is $n$-soft by \cite{FC,3.4.7}, we complete the proof
of \cref{7.5.1} as soon as the strong $n$-universality of $f$ will
be established. For this purpose, pick any $n$-dimensional Polish
space $Q$, any function $\mps{\varepsilon}{X}{(0,1)}$ assessing
closeness of maps, and any map $\mps {\varphi}QX$. Let us
construct a closed embedding $\varphi'\colon Q\hookrightarrow X$
which is $\varepsilon$-close to $\varphi$.

 Note that the space $X=\varprojlim \lbrace X_t,\eta_t\rbrace$
naturally lies in $\prod\{X_i\mid i\ge 1\}$ and
$\varphi=(\varphi_1,\varphi_2,\dots)$, where $\varphi_i$ is a map
of $Q$ into $X_i$. It is clear that
$\eta_t\circ\varphi_{t+1}=\varphi_t$, for all  $t\ge 1$.

The open cylinder $\U\subset\prod\{X_i\mid i\ge 1\}$ with the base
$\V\subset\prod\{X_i\mid 1\le i\le n\}$  and generators $\{a\times
\matop{\prod}{i>n}{}X_i\mid a\in\V\}$, being intersected with $X$,
generates the corresponding structure in $X$: the open cylinder
$\U_X\subset X$, the base $\V_X\subset \prod\{X_i\mid 1\le i\le
n\}$ and the family of generators. From the definition of
$\varprojlim \lbrace X_t\rbrace$ it easily follows that there
exists a maximal subset $\tilde{\V}\subset\prod\{X_i\mid 1\le i\le
n\}$, the intersection of which with $X$ equals the chosen base
$\V_X$:
 \item\item{$(1)$} $\tilde{\V}=X_1\times\dots
\times X_{t-1}\times \W_X$, where the set $\W_X\subset X_n$ is
open (we further identify the base $\V_X$ with this set $\W_X$).

It is easy to establish the existence of the increasing sequence
$\U_X(1)\Subset \U_X(2)\Subset\dots\subset X$ of open cylindrical
sets for which
 \item\item{$(2)$} $\cup \U_X(i)=X$;
  \item\item{$(3)$} $\U_X(2i-1)\Subset \U_X(2i)$ have
the bases $\W_X(2i-1)\Subset \W_X(2i)\subset X_{n_i}$ for all
$i\ge 1$ (we can assume without loss of generality that  $n_i=i$);
and
 \item\item{$(4)$} The map $\varepsilon$ {\it has small
oscillation on generators $\U_X(i)$}, i.e.  for any $x$ and $x'$
from one generator we have
$|\varepsilon(x)-\varepsilon(x')|<\frac{\varepsilon(x)}{10}$.

Let $\xi_i=\mps{\chi_i\circ\varphi_{i+1}}{Q}{W_i}$, where $\mps
{\chi_{i}}{X_{i+1}}{W_{i}=(Z_{i+1})_{\sigma _{i}}\times _{f_{i}}
X_{i}}$ is a characteristic map of $({\cal D}_t)$;
$A_{2i}\mean\xi_i\ii\big(Z_{i+1}\times \W_X(2i)\big)\subset Q$ and
$A_{2i-1}\mean\xi_i\ii\big(Z_{i+1}\times \W_Y(2i-1)\big)\subset
Q$. Since $\sigma'_i\circ\chi_i=\eta_i$, then
\item\item{$(5)$} $A_1\Subset A_2\Subset\dots\subset Q$
and $\cup A_i=Q$.

Fix a refining sequence of open covers $\omega_i\in\cov Q$. By
virtue of $n$-softness and $n$-completeness of the characteristic
map  $\mps {\chi_{1}}{X_{2}}{W_{1}}$, there  exists a map
$\mps{\varphi_2'}{Q}{X_2}$ such that $\chi_1\circ\varphi_2'$
equals $\xi_1=\chi_1\circ\varphi_2$, and moreover,
  \item\item{$(\dag)_1$} $\varphi_2'\restriction_{\op{Cl}A_1}$ is
an $\omega_1$-map; and
  \item\item{$(\dag)_2$}
$\varphi_2'=\varphi_2$ outside $A_2\subset Q$.

 By the same reason, there exists a
map $\mps{\varphi_3'}{Q}{X_3}$ such that $\chi_2\circ\varphi_3'$
equals $\xi_2=\chi_2\circ\varphi_3$, and moreover,
 \item\item{$(\dag)_3$} $\varphi_3'\restriction_{\op{Cl}A_3}$
is an $\omega_2$-map; and
 \item\item{$(\dag)_4$} $\varphi_3'=\varphi_3$ outside $A_4\subset
Q$.

 It should now be clear to the reader how to continue
these constructions, a result of which are the maps
$\{\varphi_i'\mid i\ge 1 \}$ (for definiteness sake we suppose
$\varphi'_1=\varphi_1$). Since
$\eta_t\circ\varphi_{t+1}'=\varphi_t'$ for all  $t\ge 1$, we have
that $\varphi'\mean(\varphi_1',\varphi_2',\varphi_3',\dots)$ is a
map passing $Q$ into $X$. It is clear that
 \item\item{$(6)$} $f\circ \varphi'=f\circ\varphi$.

Let $q\in A_{2m}\setminus A_{2m-2}$. For $1\le l\le m-1$ it
follows from $(\dag)_{2l}$  that
 \item\item{$(7)$} $\varphi_i(q)=\varphi_i'(q)$ for all  $i\le m$.

Since the oscillation of $\varepsilon$ on generators of the
cylinder $U_X(2m-2)$ is small,  $\varphi(q)$ and $\varphi'(q)$ are
$\varepsilon(\varphi(q))$-close, i.e.
$\varphi'\matop{\sim}{}{\varepsilon}\varphi$. The straightforward
check using "odd" properties $\{(\dag)_{2i-1}\}$ permits us to
assert that $\varphi'$ is a closed embedding of $Q$ into $X$.
 \enddem
 Now we track the $\op{UV}^{n-1}$-division property by a passage
to inverse limit of spectra. The following auxiliary assertion
permits us to make further inductive step.
 \proclaim{Proposition \clabel{7.5_}} Let
the commutative diagram
 $({\cal E}_t)$
 $$\CD K_{t+1} @>\theta _t>> K_t \\
 @Vh_{t+1}VV  @VV h_tV \\
Z_{t+1} @>\sigma _t>> Z_t
\endCD \tag{${\cal E}_t$} $$
(more precisely, its characteristic map) be a
$\op{UV}^{n-1}$--divider of an $n$-soft  $n$-complete map, the map
$\sigma_t$ be  $n$-soft, and the map $\mps{h_t}{K_t}{Z_t}$ be a
$\op{UV}^{n-1}$-divider of $n$-soft map $\mps{f_t}{X_t}{Z_t}$,
where $X_t\hookrightarrow K_t\in\op{UV}^{n-1}$. Then $h_{t+1}$ is
a $\op{UV}^{n-1}$-divider of an $n$-soft map
$\mps{f_{t+1}}{X_{t+1}}{Z_{t+1}}$ where $X_{t+1}\hookrightarrow
K_{t+1}\in\op{UV}^{n-1}$, moreover $\theta _t(X_{t+1})\subset
X_t$, and the commutative diagram  $({\cal D}_t)$ in which
$\eta_t\mean\theta _t\restriction_{X_{t+1}}$, is $n$-soft and
$n$-complete.
\endproclaim
 \demo{Proof of \cref{7.5_}}
Let, for definiteness sake, the characteristic map $\mps
{\chi_{t}}{K_{t+1}}{W_{t}=(Z_{t+1})_{\sigma _{t}}\times _{h_{t}}
K_{t}}$ be a $\op{UV}^{n-1}$-divider of $n$-soft and $n$-complete
map $\mps{r}{K'_{t+1}}{W_t}$ where $K'_{t+1}\hookrightarrow
K_{t+1}\in\op{UV}^{n-1}$. We will use the notation for parallel
projection: $\sigma'_t\parallel\sigma_t$ and $h_t'\parallel h_t$.

 From $n$-softness of $\sigma_t$ it follows  that
$\tilde{W}_t\mean(\sigma'_t)\ii(X_t)\hookrightarrow
W_t\in\op{UV}^{n-1}$. From here and from Proposition \cref{Def6}
(applied to the class  $\cal P$ of  all  $n$-soft strong
$n$-universal maps of Polish $\ANE(n)$-spaces) it easily follows
that
  \item\item{$(a)$} $X_{t+1}\mean \theta_t\ii(X_t)
\cap K_{t+1}'=\chi_t\ii(\tilde{W}_t)\cap K_{t+1}' \hookrightarrow
K_{t+1}\in\op{UV}^{n-1}$; and
  \item\item{$(b)$} $h_t'$ is a $\op{UV}^{n-1}$-divider of
$h_t'\restriction_{\tilde{W}_t}$.

It follows  by Proposition \cref{Def4} on composition for $\cal P$
that  the composition $h_{t+1}=h_t'\circ\chi_t$ is an
$\op{UV}^{n-1}$-divider of $\mps{f_{t+1}\mean
h_{t+1}\restriction_{X_{t+1}}}{X_{t+1}}{Z_{t+1}}$.
 \enddem
 From Propositions \cref{7.5.1} and \cref{7.5_} one can derive the
basic technical result, the further application of which for
Dranishnikov's resolution permits us to represent it as an
$\op{UV}^{n-1}$-divider of the corresponding Chigogidze's
resolution.
 \proclaim{Theorem \clabel{7.5.1+}} Let $\mps h{K}{Z}$ be a
map of inverse limits of spectra $K\mean \varprojlim \lbrace
K_t,\theta _t\rbrace$ and $Z\mean \varprojlim \lbrace Z_t,\sigma
_t\rbrace$, generated by commutative diagrams $({\cal E}_t)$,
which are $n$-conservatively  soft for all $t\ge 1$. Let also  for
all $t$,
 \item\item{$(c)$} The map $\sigma _{t}$ be $n$-soft;
 \item\item{$(d)$} The diagram $({\cal
E}_t)$ be a $\op{UV}^{n-1}$-divider of $n$-soft $n$-complete map;
and
\item\item{$(e)$} The map $\mps{h_1}{K_1}{Z_1\in\ANE(n)}$ be
$n$-conservatively soft and a $\op{UV}^{n-1}$-divider of $n$-soft
map $\mps{f_1}{X_1}{Z_1}$, where $X_1\hookrightarrow
K_1\in\op{UV}^{n-1}$.

\noindent Then the map $\mps f{X}{Z}$ of inverse limits of spectra
$X\mean \varprojlim \lbrace X_t,\eta_t\rbrace$ and $Z$, generated
by commutative diagrams $({\cal D}_t)$, $t\ge 1$, from Proposition
\cref{7.5_} \footnote{More precisely, the commutative diagram
$({\cal E}_t)$ and the map $f_t$ (from $({\cal E}_{t-1})$ for
$t>1$) generates, by \cref{7.5_}, the commutative diagram $({\cal
D}_{t}),t=1,2,\dots$.}, satisfies the following properties
 \item\item{$(f)$} $f$ is $n$-soft strongly $n$-universal with
respect to Polish spaces, and $X$, $Z$ are Polish
$\op{ANE}(n)$-spaces;
 \item\item{$(g)$} $h$ is $n$-conservatively  soft; and
 \item\item{$(h)$} $h$ is a $\op{UV}^{n-1}$-divider of $f$.
 \endproclaim
 \demo{Proof } By Proposition
\cref{7.5.1}, the  map $\mps f{X}{Z}$ is $n$-soft strongly
$n$-universal with respect to  Polish spaces. From $n$-softness of
maps $\sigma _{t}$ and $Z_1\in\ANE(n)$ it follows that
$Z\in\ANE(n)$, and hence $X$ is $\op{ANE}(n)$. The property $(g)$
follows from \cite{FC,3.4.7}.

Since inverse spectrum  $\{X_t,\eta_t\}$ consists of $n$-soft
projections, and all embeddings $X_t\hookrightarrow K_t$, by
\cref{7.5_}, are $\op{UV}^{n-1}$, it follows that
$X\hookrightarrow K\in \op{UV}^{n-1}$ which proves  $(h)$.
 \enddem


    \head\chap Proof of Theorems  \cref{1+}, \cref{1} and \cref{2}
\endhead
 To construct Dranishnikov's resolution, take a
cube $R$ of sufficiently high dimension, and represent the Hilbert
cube $R\times Q$ as the inverse limit of the spectrum
$\matop{\lim} {\longleftarrow}{}\{ Z_t\mean R\times
I^t,\sigma_t\}$ where $\mps{\sigma_t}{Z_{t+1}}{Z_t}$ is the
projection along the last factor $I$. The goal is to construct
consecutively the inverse spectrum $\{K_t,\theta_t\}$ consisting
of polyhedra whose limit is $K=\mu^k$, and the morphism
$\{\mps{h_t}{K_t}{Z_t}\}$ of these spectra which will generate
Dranishnikov's resolution $\mps {h=d_k}{\mu^k}{Z=R\times Q\cong
Q}$.

 For $t=1$ we set $K_t=R\times I^1=Z_1$ and $h_t=\Id$. Suppose that
for some  $t> 1$ there is a map $\mps{h_t}{K_t}{Z_t\mean R\times
I^t}$. The cornerstone in the construction of Dranishnikov's
resolution and the proof of Theorem \cref{1+} consists in
producing of the commutative diagram,
 $$ \CD
 K_{t+1}@>\theta _t>> K_t\\
 @V h_{t+1}VV  @VV h_tV \\
Z_{t+1} @>\sigma _t>> Z_t
\endCD \tag{${\cal E}_t$} $$
the characteristic map of which is
 \item\item{$(1)$} $k$-conservatively  soft;
 \item\item{$(2)$} A $\op{UV}^{k-1}$-divider of
$k$-complete $k$-soft map; and such that
 \item\item{$(3)$} $K_{t+1}$ admits
an $\omega$-map into  $k$-dimensional polyhedron with arbitrarily
fine cover $\omega$.

  To make sure that this the case, let us consider
the compactum  $K\mean\matop{\lim} {\longleftarrow}{}\{
K_t,\theta_t\}$ and the limit map $\mps {h=d_k}{K}{R\times Q}$,
generated by commutative diagrams $({\cal E}_t),t\ge 1$. It was
established in \cite{Dr} that  $K$ is a strong $k$-universal
$\AE(k)$-compactum of dimension  $k$, and $d_k\ii$ preserves
$\AE(k)$-spaces. From Bestvina's criterion \cite{Be} it follows
that  $K$ is homeomorphic to the Menger compactum $\mu^k$. By
Proposition \cref{7.5.1}, $h=d_k$ is $k$-conservatively  soft.

We point out that, by virtue of $(1)-(3)$, Theorem \cref{7.5.1+}
is applied. As a result we get that
 \item\item{$(i)$} The map $d_k$ is an
$\op{UV}^{k-1}$-divider of the map $\mps {c_k}{X}{Z=R\times Q}$
which is $k$-soft strong  $k$-universal with respect to Polish
spaces.

 We leave the  proof of the following strengthening of $(i)$ to
the reader.
   \proclaim{Proposition \clabel{sug+1}} Given
$\ANE(k)$-space $A\subset R$, $\mps
{c_k\restriction}{c_k\ii(A\times Q)}{A\times Q}$ is a $k$-soft
strong  $k$-universal with respect to Polish spaces, and $\mps
{d_k\restriction}{d_k\ii(A\times Q)}{A\times Q}$  is a
$\op{UV}^{k-1}$-divider of $c_k$.
 \endproclaim
 Now, it easily follows  that  $X\subset K$ is a strong
$k$-universal Polish $\AE(k)$-space of dimension $k=\dim K$. By
Characterization Theorem \cref{sug1} for N\"{o}beling space, it
follows that $X\cong \nu^k$, which {\sl completes the proof of
Theorem \cref{1+}}. The evident application of $(i)$,
\cref{Def5}$(a)$ and \cref{Def55}$(b)$ {\sl proves Theorems
\cref{1} and \cref{2}}.

 \medskip
 Now we show that  the proof of the plan
$(1)-(3)$ is easily reduced to Theorem \cref{1+!}. For this, we
consider the polyhedron $P\mean K_t\times I$ simultaneously with
an arbitrarily  fine triangulation $L$. It is clear that $P$ is a
fiberwise product of $Z_{t+1}=R\times I^{t+1}$ and $K_t$ with
respect to $\sigma_t$ and $h_t$. By assumption, there exists a
polyhedron $D$ and maps $\mps{p}{D}{P}$ and $\mps{q}{D}{P^{(k)}}$
such that \cref{1+!}$(1)-(3)$ hold. It is easily seen that
 \item\item{$(4)$} The projection $\pi$ of the compactum
$K_{t+1}\mean D\times T$ onto $D$ along the cube $T$ of dimension
$t\ge 2k+1$ is a $\op{UV}^{k-1}$-divider of $k$-soft $k$-complete
projection $D\times N^t_{k}$ onto $D$ along standard
$k$-dimensional  N\"{o}beling space  $N^t_{k}\hookrightarrow
T\in\op{UV}^{k-1}$.

 Complete the definition of the diagram  $({\cal E}_t)$ as
$\mps{\theta_t\mean\sigma_t'\circ p\circ\pi}{K_{t+1}}{K_t}$ and
$\mps{h_{t+1}\mean h_t'\circ p\circ\pi}{K_{t+1}}{Z_{t+1}}$, where
$\sigma_t'\parallel\sigma_t$ is a projection along $I$ and
$h_t'\parallel h_t$. It follows from \cref{1+!}$(2)$ and $(4)$
that the characteristic map of $(\tilde{{\cal E}}_t)$ -- the map
$\mps{p\circ\pi}{K_{t+1}}{K_t}$, is a $\op{UV}^{k-1}$-divider of a
$k$-complete $k$-soft map. It follows from \cref{1+!}$(3)$ that if
the triangulation $L$ is sufficiently fine, then the composition
$\mps{q\circ\pi}{K_{t+1}}{P^{(k)}}$ satisfies $(3)$.

   \medskip
\head\chap Multivalued retraction of a ball onto its boundary
\endhead
 In the next two sections we outline (after
\cite{Dr}) the construction of a multivalued retraction of the
ball onto its boundary (going back to I. M. Kozlowski) and
multivalued retraction of a polyhedron onto its $k$-dimensional
skeleton. We prove in passing Theorem \cref{1+!}.

Let $\partial B^{n+1}$ be the boundary of unit ball $B^{n+1},n\ge
1$. By $B_y^{n+1}$, $y\in \partial B^{n+1}$, we denote the ball of
radius $3/4$, tangent to the sphere $\partial B^{n+1}$ in $y$. It
is evident that the multivalued mappings $\Q_{n+1}\colon
\partial B^{n+1}\rightsquigarrow B^{n+1},
\Q_{n+1}(y)=B_y^{n+1}$, and $\P_{n+1}\colon
B^{n+1}\rightsquigarrow
\partial B^{n+1},\P_{n+1}(x)=
\{y\in \partial B^{n+1}\mid B_y^{n+1}\ni x\}$, are inverse each to
other. Since
\item\item{$(1)$} The restriction
$\P_{n+1}$ on $\partial B^{n+1}$ is the identity, i.e.
$\P_{n+1}(x)=x$, for every  $x\in\partial B^{n+1}$,

\noindent $\P_{n+1}$ is a multivalued retraction of the ball onto
its boundary. It is  the base of the construction of
Dranishnikov's resolution. we list several  rather easy properties
of $\P_{n+1}$ which we shall need later on:
 \proclaim{Lemma \clabel{0.5} }
\item\item{$(2)$} $\{x\in B^{n+1}\mid \P_{n+1}(x)=
\partial B^{n+1}\}=\frac{1}{2}
\cdot B^{n+1}$;
\item\item{$(3)$}  $\P_{n+1}(x)\subsetneqq \P_{n+1}(a\cdot x)$
for all  $x\in B^{n+1}\setminus\frac{1}{2}\cdot B^{n+1}$ and
$a<1$; and
\item\item{$(4)$} $\P_{n+1}(a\cdot y)\not\ni(-y)$,
for all  $y\in \partial B^{n+1}$ and $\frac{1}{2}<a\le1$.
 \endproclaim
 Next, consider the graph  $D_{n+1}\mean\{(y,x)\mid
x\in B_y^{n+1}\}\subset \partial B^{n+1}\times B^{n+1}$ of the map
$\Q_{n+1}$. It is clear that $D_{n+1}$ and the graph of $\P_{n+1}$
are symmetric with respect to the permutation of $x$- and
$y$-coordinates. By $(1)$, $\partial B^{n+1}$ is naturally
contained in $D_{n+1}$. Concerning the natural projections
$\mps{p_{n+1}}{D_{n+1}}{B^{n+1}}$ and
$\mps{q_{n+1}}{D_{n+1}}{\partial B^{n+1}}$ of the graph $D_{n+1}$
onto its factors, the following is  known \cite{ARS}:
  \proclaim{Proposition \clabel{1} } $p_{n+1}$ is
$n$-conservatively soft,  and $q_{n+1}$ is a soft retraction.
 \endproclaim
 Moreover, since  $q_{n+1}\ii(y)=B_y^{n+1}$ and $\{y\}\subset_Z
B_y^{n+1}$ $y\in \partial B^{n+1}$, it follows that
 \item\item{$(5)$} $\partial B^{n+1}\subset
D_{n+1}$ is a fiberwise $Z$-set with respect to $q_{n+1}$, i.e.
for every  partial map $\pam{Z}{A}{\varphi}{D_{n+1}}$ which is the
local lift of $\mps{\psi}{Z}{\partial B^{n+1}}$, there exists an
extension $\mps{\hat\varphi}{Z}{D_{n+1}}$ of the map $\varphi$,
which is a global lift of $\psi$, such that
$\hat\varphi(Z\setminus A)$ does not intersect $\partial B^{n+1}$.

 \smallskip
 We conclude this section by studing
the $\op{UV}^{n}$-division of $p_{n+1}$. By $C_{n+1}\subset
D_{n+1}$ we denote an open subset $D_{n+1}\setminus T^n$, where
$T^n\mean \{(y,-\frac{1}{2}\cdot y)\mid y\in \partial B^{n+1}\}$.
 \proclaim{Proposition \clabel{2}}
The map $\mps{p_{n+1}}{D_{n+1}}{B^{n+1}}$  is a
$\op{UV}^{n}$-divider of the $n$-soft map $\break
\mps{p_{n+1}\restriction_{C_{n+1}}}{C_{n+1}}{B^{n+1}}$.
 \endproclaim
 \demo{Proof of \cref{2}} The fact that  the map
$p_{n+1}\restriction_{C_{n+1}}$ of complete spaces is $n$-soft
follows  from $\{(p_{n+1})\ii (x)\cap C_{n+1}\mid x\in
B^{n+1}\}\in\op{equi}$-$\op{LC}^{n-1}$.

 The homotopy
$\mps{h_t}{B^{n+1}}{B^{n+1}},t\in I$, given by $h_t(x)=\big(1-t(1-
\parallel x\parallel)\big)\cdot x$ is called {\it radial}.
It joins $\op{Id}$ with $h_1=\parallel x\parallel\cdot x$.
Consider also  the continuous homotopy
$\mps{H_t}{D_{n+1}}{\partial B^{n+1}\times B^{n+1}},0\le t\le 1$,
given by $H(y,x)=(y,h_t(x))$. The  property
$C_{n+1}=D_{n+1}\setminus T^n\hookrightarrow
D_{n+1}\in\op{UV}^{n-1}$, which heretofore remained unproved,
easily follows from the  assertion given below. $\square$
 \proclaim{Lemma \clabel{3} } For all  $t\ge 0$,
we have $H_t(D_{n+1})\subset D_{n+1}$, and also
 \item\item{$(a)$} $H_0=\Id$ and $H_t(D_{n+1})\cap T^n=\emptyset$
for all  $t>0$.
 \endproclaim
 {\it Proof of \cref{3}.}
Since  $h_t(x)=a\cdot x$, where  $a\le 1$, it follows  by
\cref{0.5}$(iii)$ that $\P_{n+1}(x)\subset \P_{n+1}(h_t(x))$.
Hence, if $y\in\P_{n+1}(x)$, then $y\in\P_{n+1}(h_t(x))$, i.e.
$H_t(D_{n+1})\subset D_{n+1}$.

 Suppose that  $H_t(y_0,x_0)\in T^n$, for some
point $(y_0,x_0)\in D_{n+1}$,  i.e.  $h_t(x_0) =-\frac{1}{2}\cdot
y_0$, where $y_0\in \partial B^{n+1}$. Since  $h_t(x_0)=b\cdot
x_0$ where $b<1$, $x_0=\alpha\cdot(-y_0)$ for
$\alpha>\frac{1}{2}$. In view of $(4)$ we have
$\P_{n+1}(x_0)\not\ni-(-y_0)=y_0$, i.e.  $(y_0,x_0)\not\in
D_{n+1}$, a contradiction. \enddem
 \smallskip
 Fix a point (center) $O$ of the relative interior
$\op{rint}\Delta^{n+1}$ of a simplex. Then the dilation with
center $O$ generates a multiplication $a\cdot x$ for
$x\in\Delta^{n+1}$ and $0\le a\le 1$. By the {\it antipode to}
$y\in\partial\Delta^{n+1}$ we understand the intersection of the
ray $[y,O)$ with $\partial\Delta^{n+1}$. Then the multiplication
$a\cdot x$  can be extended on all $x\in\Delta^{n+1}$ and $-1\le
a\le 1$. If $x=a\cdot y$, where $y\in\partial\Delta^{n+1}$ and
$0\le a\le 1$, then $a$ is called the {\it norm} $\parallel
x\parallel$ of $x$.

Let  $\mps{\theta}{\Delta^{n+1}}{B^{n+1}}$ be a {\it radial
homeomorphism}, i.e.   $\theta(a\cdot x)=a\cdot \theta(x)$, for
all $x\in\Delta^{n+1}$ and $-1\le a\le 1$. The conjugacy operation
with respect to  homeomorphism  $\theta$ transforms all early
obtained constructions for pair  $(B^{n+1},\partial B^{n+1})$ into
constructions for pair $(\Delta^{n+1}, \partial \Delta^{n+1})$. In
particular, the radial homotopy $\mps{h_t}{B^{n+1}}{B^{n+1}}$ (see
\cref{2})  passes to {\it the radial homotopy}
$\mps{\theta\ii\circ h_t\circ \theta}{\partial
\Delta^{n+1}}{\partial \Delta^{n+1}}$ which we continue to denote
by $h_t$. Since all results obtained earlier are valid also for
simplexes, in the case of simplex we will use the previous
notations for the corresponding spaces and maps.


\head\chap Multivalued retraction of a polyhedron onto its
skeleton. Proof of Theorem \cref{1+!}
\endhead
 Let $P$ be a compact polyhedron of dimension $m$ given
with some triangulation $L$. Represent the $(n+1)$-dimensional
skeleton $P^{(n+1)},n\ge k$, as $\cup\{\Delta^{n+1}_i\mid i\ge
1\}$. By previous section, the following objects sre defined for
every  $i$: the multivalued map $\P_{n+1}(i)\colon
\Delta_i^{n+1}\rightsquigarrow
\partial\Delta_i^{n+1}$, the graph  $D_{n+1}(i)\subset
\partial\Delta_i^{n+1}
\times \Delta_i^{n+1}$ of the mapping
$\Q_{n+1}(i)=(\P_{n+1}(i))\ii$ and the natural projections
$p_{n+1}(i)\colon D_{n+1}(i)\rightarrow\Delta^{n+1}_i$ and
$\mps{q_{n+1}(i)}{D_{n+1}(i)}{\partial\Delta^{n+1}_i}$ of
$D_{n+1}(i)$ onto factors.

 Since $\P_{n+1}(i)\colon \Delta_i^{n+1}\rightsquigarrow
\partial\Delta_i^{n+1}$ and $\P_{n+1}(j)\colon \Delta_j^{n+1}
\rightsquigarrow
\partial\Delta_j^{n+1}$ for $\Delta_i^{n+1}\cap
\Delta_j^{n+1}\not=\emptyset$ agree on the common domain (where
they are identical), we have that

 \centerline{$D^n_{n+1}\mean \{(a_n,a_{n+1})
\in\partial\Delta_i^{n+1}\times \Delta_i^{n+1}\mid
(a_n,a_{n+1})\in D_{n+1}(i)\}$}

\noindent is contained in a natural manner in the  union of the
boundaries of  all simplexes $\Delta_i^{n+1}$. Also, the natural
projections $\mps{p^n_{n+1}}{D^n_{n+1}}{P^{(n+1)}}$ and
$\mps{q^n_{n+1}}{D^n_{n+1}}{P^{(n)}}$ are correctly defined. The
following is true:
 \item\item{$(i)$} $(a_n,a_{n+1})\in D^n_{n+1}$ and $a_{n+1}\in
P^{(n)}$ imply $a_n=a_{n+1}$.

\noindent  Hence $P^{(n)}$ is naturally contained in $D^n_{n+1}$,
and $q^n_{n+1}$ is a retraction. It was proved in \cite{ARS} that
$p^n_{n+1}$ is again $n$-conservatively soft map.

The homotopy $\mps{h_t}{P^{(n+1)}}{P^{(n+1)}}$ is called {\it
radial}, if its restriction on each simplex $\Delta_i^{n+1}\subset
P^{(n+1)}$ is a radial homotopy. It is clear that $h_t$ is
identity on  $P^{(n)}$. Let $C_{n+1}(i)\subset D_{n+1}(i)$ be an
open subset taken from \cref{2}, and $C^n_{n+1}\mean\cup
\{C_{n+1}(i)\mid i\}$ an open subset of $D^n_{n+1}$. The proof of
the fact that
 \item\item{$(ii)$} The map
$\mps{p^n_{n+1}}{D^n_{n+1}}{P^{(n+1)}}$  is an
$\op{UV}^{n-1}$-divider of the $n$-soft map $\break
\mps{p^n_{n+1}\restriction_{C^n_{n+1}}}{C^n_{n+1}}{P^{(n+1)}}$

\noindent is performed analogously to Proposition \cref{2} with
the help of the fact given below.
 \proclaim{Proposition \clabel{ResDr2}} The continuous
homotopy $\mps{H_t}{D^n_{n+1}}{D^n_{n+1}},0\le t\le 1$, given by
$H(y,x)
 \break=(y,h_t(x))$ transforms $D^n_{n+1}$ into $C^n_{n+1}$ for
every $t>0$.
 \endproclaim
Since $q^n_{n+1}$ is not open,  $q^n_{n+1}$ fails to possess the
softness property which is inherent in $q_{n+1}(i)$. But
nevertheless a weak softness property of $q^n_{n+1}$  can be
detected.
 \proclaim{Proposition \clabel{Res1} } The projection
$q^n_{n+1}$ is  a synchronized  Hurewich fibration.
 \endproclaim
 We introduce the necessary notions.
 \definition{Definition  \clabel{Def2}}
The homotopy $\mps {\varphi}{X\times I}{P^{(s)}}$ is called {\it
synchronized } if
 \item\item{$(1)$} $\varphi\ii(\Delta)=(\varphi_0)
\ii(\Delta)\times I$, for every simplex $\Delta\subset P^{(s)}$.
 \enddefinition
 In other words, $(1)$ means that, if $\varphi_{0}(x)\in\Delta$,
then $\varphi_t(x)\in\Delta$  for all  $t\in I$.
 \definition{Definition  \clabel{Def3} }
The map $\mps f{D^n_{n+1}}{P^{(n)}}$ is called {\it synchronized
Hurewich fibration } if for every synchronized  homotopy $\mps
{\varphi}{X\times I}{P^{(n)}}$ and for every partial lift  $\mps
{\theta_0}{X\times \{0\}}{D^n_{n+1}}$ of map $\varphi_0$ with
respect to  $f$ there exists a homotopy $\mps {\theta}{X\times
I}{D^n_{n+1}}$ lifting the homotopy $\varphi$ such that $\mps
{p_{n+1}^n\circ \theta}{X\times I}{P^{(n+1)}}$ is a synchronized
homotopy.
 \enddefinition
 \demo{Proof of \cref{Res1}}
 Let $\mps {\varphi}{X\times I}{P^{(n)}}$ be a synchronized
homotopy and $\mps {\psi_0}{X\times \{0\}}{P^{(n+1)}}$ a map. For
the proof it is sufficient to establish that if the map
$\theta_0\mean(\varphi_0,\psi_0)$ transforms $X$ into $D^n_{n+1}$,
then there exists a synchronized  homotopy $\mps {\psi}{X\times
I}{P^{(n+1)}}$ extending $\psi_0$ such that
$\theta\mean(\varphi,\psi)$ is a homotopy of $X$ into $D^n_{n+1}$.

Let $P^{(n+1)}= \{\Delta_i^{(n+1)}\mid i\ge 1\}$. Consider the
following subsets of $X$: $X_0\mean (\psi_0) \ii(P^{(n)})$ and
$X_i\mean(\psi_0) \ii(\Delta_i^{(n+1)})$. It is clear that $X=\cup
X_i$ and
 \item\item{$(iii)$} $X_i\setminus X_0\subset\op{Int}X_i$
for all  $i\ge 1$.

It follows from $(i)$  that  $\varphi_0=\psi_0$ on $X_0$. As
$\psi_0(X_i)\subset \Delta_i^{(n+1)}$ and
$\op{Im}(\theta_0)\subset D^n_{n+1}$, then
$\varphi_0(X_i)\subset\partial\Delta_i^{(n+1)}$. Since $\varphi$
is the synchronized  homotopy, we have
$\varphi_t(X_i)\subset\partial\Delta_i^{(n+1)}$, for every $t\in
I$.

 Given  $i\ge 1$, consider the following commutative diagram,
  {$$

\def\p#1,#2 {\llap{$\ssize p_{#1}^{#2}$}}
\matrix
 &D_{n+1}(i)&
 \matop{\longrightarrow}{}{q_{(n+1)}(i) } &
\partial\Delta_i^{(n+1)}\cr
 &\uparrow\sigma_i&&\uparrow\varphi\cr
 &A_i&\matop{\hookrightarrow}{}{} &X_i\times I\endmatrix
 $$}
\noindent in which $A_i\mean \big((X_i\cap X_0)\times
I\big)\cup\big(X_i\times\{0\}\big)$, $\sigma_i=\theta_0$ on
$X_i\times\{0\}$ and $\sigma_i=(\varphi,\varphi)$ on $(X_i\cap
X_0)\times I$. Since $q_{n+1}(i)$ is soft, there exists an
extension $\mps{\theta_i}{X_i\times I}{D_{n+1}(i)}$ of $\sigma_i$
such that  $q_{n+1}(i)\circ \theta_i=
\varphi\restriction_{X_i\times I}$.

 By Proposition \cref{1},
$\partial\Delta_i^{(n+1)}\subset D_{n+1}(i)$ is a fiberwise
$Z$-set with respect to  $q_{n+1}(i)$. Then $\theta_i$ can be
chosen in a such manner that
 \item\item{$(iv)$}  $\theta_i(X_i\times I\setminus A_i)$
is contained in  $D_{n+1}(i)\setminus\partial\Delta_i^{(n+1)}$
(i.e. $p^n_{n+1}\circ \theta_i(X_i\times I\setminus
A_i)\subset\op{rint} \Delta_i^{(n+1)}$).

 The desired homotopy $\mps {\theta}{X\times I}{D^n_{n+1}}$
equals $\theta_i$ on $X_i\times I$. We can check straightforwardly
with help of $(iii)$ and $(iv)$ that $\theta$ is continuous, and
$p_{n+1}^n\circ \theta$ is a synchronized  homotopy.
 \enddem


 For any $k\le s<t\le m$ consider the increasing
sequence $P^{(s)}\subset P^{(s+1)}\subset \dots\subset
P^{(t-1)}\subset P^{(t)}$ of the skeleta of the polyhedron $P,\dim
P=m$. In the general situation let us define the following
objects: $D_t^s\mean\lbrace a=(a_s,a_{s+1},\dots,a_{t-1},a_{t})\in
P^{(s)}\times \dots \times P^{(t-1)}\times P^t\mid
(a_i,a_{i+1})\in D^{i}_{i+1}, s\le i<t\rbrace \subset
P^{(s)}\times \dots \times P^{(t-1)}\times P^t$, the maps
$\mps{p_t^s}{D_t^s}{P^{(t)}}$  and $\mps{q_t^s}{D_t^s}{P^{(s)}}$
by formulas $p_t^s\mean a_{t}$ and $q_t^s\mean a_s$ respectively.
It was proved in \cite{Dr} that $D_t^s\in\ANE$.

 \proclaim{Conjecture \clabel{C} } The compactum   $D_t^s$ is a polyhedron.
 \endproclaim
 
 This plausible (and, apparently, difficult) conjecture
was formulated (but not proved) in \cite {Dr, p.124}. It seems
likely  that a reason of this lies in the fact that it was easier
to make use of already known results of $Q$-manifolds theory. We
also do not want to spend effort on the proof of this conjecture
as the basic result of the paper does not depend on its validity
(in the case of conjecture failure, one must draw on the Edwards
Theorem and Chapman Theorem as made in \cite{Dr}). But for the
simplicity of the text we do assume  that {\it $D_t^s$ is a
polyhedron}.

 Analogously to \cref{Def3} we say that  the map
$\mps{q_t^s}{D_t^s}{P^{(s)}}$ is an {\it synchronized  Hurewich
fibration}, if for every synchronized  homotopy $\mps
{\varphi}{X\times I}{P^{(s)}}$ and for every  partial lift $\mps
{\theta_0}{X\times \{0\}}{D_t^s}$ of $\varphi_0$ there exists a
homotopy  $\mps {\theta}{X\times I}{D_t^s}$ which lifts  $\varphi$
such that $p_t^s\circ \theta$ is a synchronized  homotopy.
Applying Proposition \cref{Res1} several times, it can be easily
proved that
 \item\item{$(v)$} The map  $q_t^s$ (as well as $q^n_{n+1}$)
is a synchronized  Hurewich fibration.

 \medskip
  By $C_t^s$ we denote an open subset $\{ a\mid
(a_i,a_{i+1})\in C_{i+1}^i\ \text{for all }\ i,s\le i<t\}\subset
D_t^s$ ($C_{i+1}^i\subset D_{i+1}^i$ is taken from $(ii)$). The
essential complement of \cite{ARS} where the $k$-conservative
softness of  $p_t^s$ was established is the following {\it key
result} which proves Theorem \cref{1+!}:

 \proclaim{Theorem \clabel{Res2} } The map
$\mps{p_m^k}{D_m^k}{P=P^{(m)}}$ is a $\op{UV}^{k-1}$-divider of
the $k$-soft map
$\break\mps{p_m^k\restriction_{C_m^k}}{C_m^k}{P}$.
 \endproclaim
 
 \demo{Proof of \cref{Res2}}
 Consider a closed subset

 \centerline{$F=D_m^k\setminus C_m^k=
\{a\mid (a_i,a_{i+1})\not\in C_{i+1}^i\ \text{for some}\ i,k\le
i<m\}\subset D_m^k$}

\noindent and its closed filtration $\ \ F_m\subset
F_{m+1}\subset\dots\subset F_{k+2}\subset  F_{k+1}=F\ $ where

 \centerline{$F_s=\{a\in F\mid (a_i,a_{i+1})\in
C_{i+1}^i\ \text{for all }\ i,k\le i<s-1\}$.}

 \proclaim{Lemma \clabel{Res3}} The restriction
$\mps{p_m^k\restriction_{C_m^k}}{C_m^k}{P=P^{(m)}}$ is a $k$-soft
map.
 \endproclaim
 
 \demo{Proof of \cref{Res3}} Suppose that the partial map
$Z\supset A\matop{\rightarrow}{}{\varphi}C_m^k
\matop{\rightarrow}{}{p_m^k}P,\dim Z\le k$, has an extension
$\mps{\psi}{Z}{P}$. Represent $\varphi$ in the coordinate form
$(\varphi_k,\varphi_{k+1},\dots,\varphi_m)$ where $\varphi_i$ is
the map from $A$ into $P^{(i)}$. Then
$\varphi_m=\psi\restriction_{A}$. Since
$\mps{p_{i+1}^i}{C_{i+1}^i}{P^{(i+1)}}$ is $k$-soft for all $i$,
we can construct, by inverse induction on $m,m-1,\dots,k+1$, the
maps $\psi_m=\psi,\psi_{m-1},\dots,\psi_{k+1}$ from $Z$ into
$P^{(i)}$ such that  $\varphi_i=\psi_i\restriction_{A}$ and
$(\psi_{i+1},\psi_i)\in C_{i+1}^i$ for all  $i$. In view of
$p_m^k\circ \hat\varphi=\psi$ it is clear that
$\hat\varphi\mean(\psi_k,\psi_{k+1},\dots,\psi_m)\colon
Z\rightarrow C_m^k$ is the desired extension of  $\psi$.
\enddem
 To complete the proof of
Theorem \cref{Res2} it is sufficient to show that  the map
$\mps{\varphi}{A}{D_m^k}$ of any compactum  $A$  is arbitrarily
closely approximable by a map which does not intersect $F$. Let
$\mps{h_t}{P^{(m)}}{P^{(m)}}$ be a radial homotopy, and let
$\mps{H_t}{D_m^k}{D_m^k}$ a homotopy given by $\break
H_t(a_k,a_{k+1},\dots,a_{m-1},a_m)=(a_k,a_{k+1},
\dots,a_{m-1},h_t(a_m))$. It is clear that the homotopy
$H_t\circ\varphi$ removes $A$ from $F_m$.

 Taking this remark into account, it is sufficient to show that if
$\varphi(A)\cap F_{s+1}=\emptyset, k\le s<m$, then  there exists a
map $\mps{\varphi'}{A}{D^k_m}$ arbitrarily close  to  $\varphi$
such that  $\varphi'(A)\cap F_{s}=\emptyset$. Again represent
$\varphi$ in the coordinate form
$(\varphi_k,\varphi_{k+1},\dots,\varphi_m)$, where $\varphi_i$ is
the map from $A$ into $P^{(i)}$. Let $\mps{h_t}{P^{(s)}}{P^{(s)}}$
be a radial homotopy. Then $\mps{\Psi_t^s\mean
h_t\circ\varphi_s}{A}{P^{(s)}},0\le t\le 1$, is a synchronized
homotopy. As was noted in $(v)$, the map
$\mps{q^s_{m}}{D^s_{m}}{P^{(s)}},q^s_{m}(a)=a_s$, is an
synchronized Hurewich fibration. Therefore, there exist
synchronized homotopies $\Psi_t^{s+1},\Psi_t^{s+2},\dots,\Psi_t^m$
from $A$ into $P^{(s+1)},P^{(s+2)},\dots,P^{(m)}$ such that  the
formula
$\Psi_t\mean(\Psi_t^s,\Psi_t^{s+1},\Psi_t^{s+2},\dots,\Psi_t^m)$
defines the homotopy from $A$ into $D_m^s$.

We take

 \centerline{$(\varphi_k,\varphi_{k+1}\dots,\varphi_{s-1},\Psi_t)=
(\varphi_k,\varphi_{k+1}\dots,\varphi_{s-1},
\Psi_t^s,\Psi_t^{s+1},\Psi_t^{s+2},\dots,\Psi_t^m)$}

\noindent as a homotopy $\mps{\Phi_t}{A}{D_m^k}$ removing $A$ from
$F_{s}$. Since $(\varphi_{s-1},\Psi_t^s)(A)\subset D^{s-1}_s$, for
every $t\ge 0$, we easily deduce $\Phi_t(A)\subset D_m^k$. Next,
we note that, due to Proposition \cref{ResDr2},
$(\varphi_{s-1},\Psi_t^s)$ maps $A$ into $C^{s-1}_{s}$ for every
$t>0$.  Hence, we have proved that for every  $\delta>0$,
$\varphi\matop{\sim}{}{\delta}\Phi_t$, for sufficiently small
$t>0$, and $\op{Im}(\Phi_t)\cap F_{s}=\emptyset$.
 \enddem
 
 \medskip
\head\chap Epilogue
\endhead
 Here we list a selection of applications of
results obtained in the paper, and also formulate some unsolved
problems.

 {\it Uniqueness problem of Chigogidze's resolution.}
 By the {\it $k$-dimensional
Chigogidze's resolution  over $Y\in\ANE(k)$} we understand a
$k$-soft map $\mps{f}{X}{Y}$ of $k$-dimensional space  $X$ onto
$Y$, which is strongly $k$-universal with respect to maps of
Polish spaces. One of the central problems of the N\"{o}beling
spaces theory consists in establishing of the topological
uniqueness of such a resolution \cite{FC}.

 \proclaim{Problem \clabel{8.1}} Prove  that  any two Chigogidze's resolution
$\mps{f,g}{\nu^k}{Q}$ are homeomorphic, i.e.  there exists a
homeomorphism $\mps{h}{\nu^k}{\nu^k}$ such that $f=g\circ h$.
\endproclaim

 For $k=\infty$ this problem was solved in
affirmative \cite{Tor}. The case $k=0$ was also settled (see, for
example, \cite{AGS}).
 \smallskip
 {\it Problem of the characterization of Dranishnikov's resolution.}
This resolution no doubt represent the analogy of the Menger
compactum in the category of maps. In analogy with compacta, the
question of its  characterization arises naturally. But prior to
doing  this, we should understand what is Dranishnikov's
resolution. In view of the results of this paper, the {\it
$k$-dimensional Dranishnikov's resolution over $Y\in\ANE(k)$} is
any proper map $\mps{f}{X}{Y}$ from a $k$-dimensional space $X$
onto $Y$ such that
 \roster
 \item $f$ is $k$-conservatively soft
strong $k$-universal with respect to compacta; and
 \item $f$ is a $\op{UV}^{k-1}$-divider of $k$-dimensional
Chigogidze's resolution over $Y$.
\endroster

From $k$-conservative softness of $f$ it follows that
 \item\item{$(3)$}  $f$ is polyhedral $k$-soft, $(k,k-2)$- and
$(k-1)$-soft;
\item\item{$(4)$} $f\ii$ preserves $\op {AE}(k)$-spaces, and
therefore, $f$ is a $\op{UV}^{k-1}$-map;
 \item\item{$(5)$} $f$ is $k$-inversible map;
 \item\item{$(6)$} The preimage  of any
conservatively closed $\op {equi}$-$\op {LC}^{k-1}$-family
$\lbrace Y_\alpha \rbrace$  is $\op {equi}$-$\op {LC}^{k-1}$.

 From strong  $k$-universality of proper map $f$ with
respect to  compacta it follows  that
 \item\item{$(7)$}  $f$ is $k$-universal
map with respect to  Polish spaces.

From the fact that  $f$ is a $\op{UV}^{k-1}$-divider of
$k$-dimensional  Chigogidze's resolution   it follows  that
 \item\item{$(8)$}  $f\ii$  preserves $Z$-sets and
strong  $k$-universal spaces.

 There is a definite hope that  the topological type
of Dranishnikov's resolution is unique.
 \proclaim{Problem \clabel{8.2}} Are any two
Dranishnikov's resolutions $\mps{f,g}{\mu^k}{Q}$ homeomorphic?
\endproclaim

 {\it Triangulation Theorem for $\nu^k$-manifolds.}
 In a natural way one can define manifolds modelled on the
N\"{o}beling space ($\nu^k$-manifolds) and on the universal
pseudoboundary ($\sigma^k$-manifolds). It follows from
Propositions \cref{sug+1} and \cref{Def6}  that  there exists a
$k$-dimensional Dranishnikov's  resolution $\mps{f}{X}{Y}$ over
any $l_2$-manifold (over any $\Sigma$-manifold where $\Sigma$ is
the universal pseudoboundary of Hilbert cube) $Y$. By virtue of
the Characterization Theorem for $\nu^k$-manifolds (see \cite{AgM}
and \cite{Ng}), it follows that $X$ is homeomorphic to
$\nu^k$-manifold. The converse fact is true that is the content of
the Triangulation Theorem:
\item\item{$(a)$} For every
$\nu^k$-manifold  $X$ there exists a $k$-dimensional
Dranishnikov's resolution  $\mps{f}{X}{Y}$ in which $Y$ is an
$l_2$-manifold.

 \noindent Since any $l_2$-manifold
is homeomorphic to an open subset of $l_2$, the central result of
Nag\"{o}rko thesis \cite{Ng} follows easily
\item\item{$(b)$} Any $\nu^k$-manifold
is homeomorphic to an open subset of $\nu^k$.

 Analogous results take place for $\sigma^k$-manifolds.

It is easy to establish the following relation between $\mu^k$-,
$\nu^k$- and $\sigma^k$-manifolds.
 \item\item{$(c)$} Any
$\nu^k$-manifold admits a $\op{UV}^{k-1}$-embedding into a
$\mu^k$-manifold; and
 \item\item{$(d)$} Any $\sigma^k$-manifold lies in a
$\nu^k$-manifold as ${\cal C}_{c(k)}$-absorber set  where ${\cal
C}_{c(k)}$ is the class of all $k$-dimensional compacta.

 The following general question remains unsolved.
 
 \proclaim{Problem  \clabel{8.3}} When is the preimage  $d_k\ii(X), \
X\subset Q$, homeomorphic to $\nu^k$-manifold.
 \endproclaim

 {\it Problem of geometrization of Dranishnikov's resolution}.
 Initially Dranishnikov's and Chigogidze's resolutions were
constructed in a unconstructive manner as the limit projections of
some countable spectrum. We can identify their domains lying in
Hilbert cube with Menger and N\"{o}beling spaces only with help of
corresponding Characterization Theorems. On the other hand, in
\cite{Ag2} Chigogidze's resolution  was constructed in a geometric
manner as the orthogonal projection of the standard N\"{o}beling
space. The fractal structure of this resolution  was thereby
revealed. It was interesting to realize Dranishnikov's resolution
also in a geometric manner. We precede the formulation of the
corresponding conjecture by the series of definitions.

{\it The standard Menger space $M^{m}_k$} and {\it geometric
pseudointerior  $I(M^{m}_k)$} can be defined with help of the
notion of $\Sigma$-product. Recall that by $\Sigma^n_k(X;Y)\subset
X^n$, where $Y\subseteq X$ and $k\le n<\infty$ we denote the union
of images $Y^{n-k}\times(X\setminus Y )^{k}\subset X^n$ permuting
by all coordinates. Given $m\ge 2k+1$,
$M^{m}_k\mean\bigcap\{\Sigma^{m}_k\left([0,1];\C_n\right)\mid 1\le
n<\infty\}\subset I^m$ and
$I(M^{m}_k)\mean\bigcap\{\Sigma^{m}_k\left([0,1];{\cal
J}_n\right)\mid 1\le n<\infty\}\subset M^{m}_k$ where ${\cal
J}_n\mean [0,3^{-n})\cup\cup\{3^{-n}\cdot (3i-1,3i+1)\mid 1\le
i<3^{n-1}\}\cup (1-3^{-n},1]\subset\C_n\mean
[0,3^{-n}]\cup\cup\{3^{-n}\cdot [3i-1,3i+1]\mid 1\le
i<3^{n-1}\}\cup [1-3^{-n},1]$. The standard Menger space $M^{m}_k$
and geometric pseudointerior  $I(M^{m}_k)$ are homeomorphic to
$\mu^k$ and $\nu^k$, respectively.

 \proclaim{Conjecture \clabel{8.4}} Let $m\ge(2k+1)+(k+1)^2$.
Is it true that  there exists an orthogonal projection
$\mps{p}{\Bbb R^m}{\Sigma}$ onto $(2k+1)$-dimensional subspace
$\Sigma$ such that $\mps{p\restriction}{M^m_k}{p(M^m_k)}$ has the
same soft properties as Dranishnikov's resolution? Is it
$\op{UV}^{k-1}$? Is it true that
$\mps{p\restriction}{I(M^m_k)}{p(I(M^m_k))}$ is a Chigogidze's
resolution? Is it true that  $p\restriction_{M^m_k}$ is a
$\op{UV}^{k-1}$-divider of $p\restriction_{I(M^m_k)}$?
 \endproclaim

 {\it Problem of the $k$-soft core.} In \cite{ARS}
the nonhomogeneity of Dranishnikov's resolution $d_k$ was
established which breaks the analogy with the case of spaces. Such
a property follows from the results on {\it soft core} and the
fact that $d_k$ is $(k-1)$- but is not $k$-soft. By the
Finite-Dimensional Michael Selection Theorem, the latter means
that the $k$-soft core

 \centerline{$\op {{\goth s}}_k(d_k)\mean\lbrace x\in\mu^k\mid
\text{collection of all fibers}\ d_k\ \text{in}\ x\ \text{is
equi-local}\ (k-1)\text{-connected}\rbrace$}

\noindent does not coincide with $\mu^k$. With the help of
additional analysis we can show that Chigogidze's resolution, a
$\op{UV}^{n-1}$-divider of Dranishnikov's resolution  $d_k$, is
contained, in fact, in the $k$-soft core $\op {{\goth s}}_k(d_k)$.
In this connection the series of questions arises.
 \proclaim{Problem \clabel{8.5}} Is it true that
$k$-soft core $\op {{\goth s}}_k(d_k)$ is homeomorphic $\nu^k$;
the restriction of $d_k$ on this set is Chigogidze's resolution
$c_k$ (and, therefore, $d_k$ is a $\op{UV}^{k-1}$-divider of this
restriction of $c_k$)?
\endproclaim
 Since for every polyhedron $P\subset Q,\dim P\le k$,
$\mps{d_k}{\mu^k}{Q}$ is $k$-inversible, there exists a section
$\mps{s}{P}{\mu^k}$. The following question is concerned with the
possibility of constructing the section  $s$ in the
equi-continuous manner, in the following sense.

 \proclaim{Problem \clabel{8.6}}
For every  $\varepsilon>0$ there exists $\delta>0$ such that for
any polyhedron  $P\subset Q,\dim P\le k$, there exists a section
$s$ of  $d_k$ such that $\diam{s(A)}<\varepsilon$, for every
$A\subset P$ with $\diam{A}<\delta$.
 \endproclaim

\newpage
\medskip


\Refs \widestnumber\key{AJNN}

 \ref
 \key Ag
 \by  S. M. Ageev
  \paper Axiomatic method of partitions in the theory of Menger and N\" obeling spaces
 \jour preprint, http://at.yorku.ca/v/a/a/a/87.htm
  \vol
  \issue
  \yr 2000
  \pages
 \endref
 
  \ref 
 \key AgM
 \by  S. M. Ageev
  \paper Axiomatic method of partitions in the theory of N\" obeling spaces, I--III
 \jour Math. Sbornik
  \vol
  \pages to appear
 \endref
 
  
 \ref
  \key Ag2
  \by S. M. Ageev
   \paper Geometric Chigogidze's resolution
  \jour submitted to Fund. Math
   \vol
   \yr
 \endref
 
 \ref
  \key Ags
  \by S. M.  Ageev, G. Gruzdev and Z. Silaeva
   \paper Characterization of 0-dimensional Chigogidze's resolution
  \jour Vestnik BGU
   \vol
   \pages to appear (in Russian)
 \endref
 
 \ref
 \key ARS
 \by  S. M. Ageev, D. Repov\v{s} and E. V. Shchepin
  \paper On the  softness of the Dranishnikov resolution 
  \jour Proc. Steklov Inst. Math.
   \vol 212
  \issue 1
  \yr 1996
  \pages 3--27
 \endref
   

 
\ref
 \key Be
 \by M.  Bestvina
  \paper  Characterizing k-dimensional universal Menger compacta
 \jour Mem. Amer. Math. Soc.
  \vol 71
  \issue 380
  \yr 1988
 \endref
 
\ref
 \key BM
 \by M.  Bestvina and J. Mogilski
  \paper  Characterizing certain incomplete infinite-dimensional  absolute retracts
 \jour Michigan Math. J.
  \vol 33
  \yr 1986
\pages  291--313
 \endref
 
 \ref 
 \key Bw
 \by  P. L. Bowers
  \paper Dense embedding of sigma-compact, nowhere locally   compact metric spaces
 \jour  Proc. Amer. Math. Soc.
  \vol 95
  \issue 1
  \yr 1985
  \pages  123--130
  \endref
  
 \ref
 \key C
 \by  A. Chigogidze
  \paper  $\op{UV}^n$-equivalence and $n$-equivalence
 \jour Topology  Appl.
  \vol 45
  \yr 1992
  \pages  283--291
 \endref


 \ref
 \key CKT
 \by  A. Chigogidze, K. Kawamura and  E. D. Tymchatyn
  \paper  N\" obeling spaces and the pseudo-interiors of Menger compacta
 \jour Topology Appl.
  \vol 68
  \yr 1996
  \pages  33--65
 \endref

\ref
 \key CZ
 \by A. Chigogidze and M. Zarichnyi
\paper Universal N\" obeling spaces and pseudo-boundaries of Euclidean spaces
 \jour Mat. Stud.
  \vol 19
  \yr 2003
 \issue 2
  \pages 193--200
 \endref
 

\ref
 \key Dr
 \by  A. N. Dranishnikov
  \paper Universal Menger compacta and universal mappings 
 \jour  Math. Sbornik
  \vol  129
  \issue 1
  \yr 1986
  \pages  121--139 (in Russian)
 \endref
  
 
 
 \ref 
 \key Eng
 \by  R. Engelking 
 \paper  General Topology
   \publaddr PWN, Warsaw
  \yr 1977
 \endref
 
 \ref 
 \key FC
 \by  V. V. Fedorchuk and A. Ch.  Chigogidze
 \paper Absolute Retracts and Infinite-dimensional Manifolds 
  \publaddr  Nauka, Moscow (in Russian)
  \yr 1992 
  \endref
  
  \ref 
  \key Hu
 \by  S. T. Hu
  \paper Theory of Retracts
  \publaddr Wayne State Univ. Press
  \yr 1965
  \endref
  
\ref
\key Ng
 \by A.   Nag\'{o}rko
  \paper Carrier and nerve theorems in the extension theory
 \jour Proc. Amer. Math. Soc.
  \vol  135 
  \yr  2007
  \issue  2
  \pages  551--558
 \endref
  
 \ref
  \key Sch
  \by E. V. Shchepin
   \paper  On homotopically regular mappings of manifolds
  \book  Geometric and Algebraic Topology,  Banach Center Publ., Vol. {\bf 18}
   \publaddr PWN, Warsaw
   \yr 1986
   \pages 139--151
 \endref
 
 \ref
 \key SB
 \by E. V. Shchepin and N. B. Brodskyi
  \paper Selections of filtered multivalued mapping
 \jour Proc. Steklov Inst. Math.
  \vol 212
  \yr 1996
  \issue 1
  \pages  209--229 
 \endref
 
 
 
 \ref
 \key {Tor}
  \by  H. Toru\' nczyk  
  \paper  Characterizing Hilbert space topology
 \jour Fund. Math.
  \vol 111
  \yr 1981
  \pages  247--262
\endref

 
 
 \ref
 \key MR
 \by J. van Mill and G. M. Reed, Eds.
  \book  Open Problems in Topology
 \publaddr North-Holland, Amsterdam
  \yr 1990
  \endref
  
 \ref
 \key Za
 \by M. M.  Zarichnyi
\paper Absorbing sets for $n$-dimensional spaces in absolutely Borel and projective classes
 \jour Math. Sbornik 
  \vol 188
  \yr 1997
 \issue 3
  \pages 435--447
 \endref
 
 \endRefs
\enddocument